\documentclass[12pt]{amsart}

\usepackage[bookmarks=true]{hyperref}
\hypersetup{colorlinks,citecolor=blue,linkcolor=blue}

\usepackage{mathptmx}
\usepackage{amsmath}
\usepackage{amssymb}
\usepackage{amsthm}
\usepackage{array}
\usepackage[mathscr]{eucal}
\usepackage[all,tips]{xy}
\usepackage{enumerate}
\usepackage{graphicx}

\newcommand{\B}[1]{\ensuremath{\mathbf{#1}}}

\newcommand{\Cal}[1]{\ensuremath{\mathcal{#1}}}

\newtheorem{Theorem}{Theorem}

\newtheorem{lemma}[subsection]{Lemma}
\newtheorem{prop}[subsection]{Proposition}

\theoremstyle{definition}

\DeclareMathOperator{\Spin}{Spin}
\DeclareMathOperator{\SO}{SO}
\DeclareMathOperator{\PO}{PO}

\newcommand\isom{\operatorname{Isom}(\mathbb{H}^n)}
\newcommand\isomp{\operatorname{Isom^+}(\mathbb{H}^n)}
\newcommand\Hy{\mathbb{H}}

\newcommand\LieH{\mathfrak{h}}

\newcommand\Z{{\mathbb Z}}
\newcommand\Q{{\mathbb Q}}
\newcommand\R{{\mathbb R}}
\newcommand\C{{\mathbb C}}
\newcommand\A{{\mathbb A}}

\newcommand\vol{{\rm Vol}_\Hy}

\newcommand\cP{{\mathcal P}}
\newcommand\Pa{{\mathrm P}}
\newcommand\bPa{{\overline{\mathrm P}}}

\newcommand\cC{{\mathcal C}}
\newcommand\D{\mathcal{D}}
\newcommand\dm{{\rm dim\;}}
\newcommand\cE{{\mathcal E}}
\newcommand\Mv{{\overline{\rm M}_v}}
\newcommand\cMv{{\overline{\mathcal M}_v}}
\newcommand\fv{{\mathfrak f}_v}

\newcommand\G{{\mathrm G}}

\newcommand\bG{{\overline{\mathrm G}}}
\newcommand\CG{{\mathrm C}}
\newcommand\gH{H}
\newcommand\bH{\overline{\gH}}

\newcommand\HH{{\mathrm H}}

\newcommand\Aut{{\rm Aut}}
\newcommand\End{{\rm End}}
\newcommand\Res{{\rm R}}
\newcommand\No{N}
\newcommand\nc{{\mathfrak n}}
\newcommand\mmax{{\mathfrak m}}
\newcommand\LL{{\mathbf L}}
\newcommand\AL{{\mathbf A}}
\newcommand\V{{\mathbf V}}
\newcommand\h{{\mathbf h}}

\newcommand\kun{\hat{k}}

\newcommand\emb[2]{{#2^{#1}}}

\newcommand\Bn{\mathrm{B}}

\newcommand\Dn{\mathrm{D}}

\begin{document}
\title{On volumes of arithmetic quotients of $\PO(n,1)^\circ$, $n$ odd}

\author{Mikhail Belolipetsky}
\thanks{Belolipetsky partially supported by EPSRC grant EP/F022662/1}
\address{
IMPA\\
Estrada Dona Castorina, 110\\
22460-320 Rio de Janeiro, Brazil}
\email{mbel@impa.br}

\author{Vincent Emery}
\thanks{Emery partially supported by SNSF, project no.
200020-121506/1 and fellowship no. PBFRP2-128067}
\address{
Department of Mathematics\\
University of Fribourg, P\'erolles\\
Chemin du Mus\'ee 23\\
CH-1700 Fribourg\\
Switzerland
}
\email{vincent.emery@gmail.com}


\subjclass[2010]{22E40 (primary); 11E57, 20G30, 51M25 (secondary)}

\begin{abstract}
We determine the minimal volume of arithmetic hyperbolic orientable $n$-dimensional orbifolds (compact and non-compact) for every odd dimension $n\ge 5$. Combined with the previously known results it solves the minimal volume problem for arithmetic hyperbolic $n$-orbifolds in all dimensions.
\end{abstract}

\maketitle

\setcounter{tocdepth}{1}

\section{Introduction}
\label{section:intro}
In this paper we study volumes of arithmetic hyperbolic orbifolds
of odd dimension. Our main results determine the arithmetic
hyperbolic orientable $n$-dimensional orbifolds (compact and
non-compact) of smallest volume  for every odd dimension $n \ge
5$. This work can be considered as a continuation  of the
research which was done in \cite{CF} for $n=3$ and in \cite{B1}
for even $n>2$.  The question of minimal volume for $n=2$ is
answered by a classical theorem of Siegel \cite{Siegel}. Together
with these results our paper concludes the study of the minimal
volume for arithmetic hyperbolic $n$-orbifolds in all dimensions.
Our Theorems~1 and 2 give a more precise version of the results of
the second author's PhD thesis~\cite{E-PhD}. We refer to the thesis
for an introduction to the methods used in this paper and for a
more detailed exposition of the proofs.

\medskip

Let $\Hy^n$ be  the hyperbolic $n$-space and let  $\isom$ be its group
of isometries.  Given a lattice $\Gamma$ in $ \PO(n,1) = \isom $, the
corresponding    quotient   space    $\Hy^n/\Gamma$   is    called   a
\emph{hyperbolic $n$-orbifold}.  When $\Gamma <  \isomp = \PO(n,1)^\circ$,
 the orbifold is orientable.
If  $\Gamma$  is an  arithmetic  lattice
defined over a number  field $k$ (see Section~\ref{arithmetic-sbgp-def}), we call
$O =  \Hy^n/\Gamma$ an \emph{arithmetic orbifold} and  we call  $k$ its
\emph{field of definition}.

\begin{Theorem}
  For each  dimension $n = 2r-1  \ge 5$, there is  a unique orientable
  compact arithmetic  hyperbolic $n$-orbifold $O^n_0$  of the smallest
  volume. It is defined over $k_0 = \Q[\sqrt{5}]$ and has hyperbolic
  volume given by the formula
\begin{eqnarray*}
   \vol(O^n_0) &=&  \frac{5^{r^2-r/2}
     \cdot 11^{r-1/2} \cdot (r-1)!}{2^{2r-1} \pi^r} \; L_{\ell_0|k_0}\!(r) \; \prod_{i=1}^{r-1} \frac{(2i -1)!^2}{(2
     \pi)^{4i}} \zeta_{k_0}(2i),
  \end{eqnarray*}
  where  $\ell_0$ is the quartic field with a defining polynomial $x^4-x^3+2x-1$.
\end{Theorem}

The proof of  Theorem 1 is conceptually similar to that in
\cite{B1}: The main ingredients are Prasad's volume formula
\cite{P}, some results and ideas of Borel-Prasad \cite{BP} and
Bruhat-Tits theory (see \cite{T2}). Essential differences in
comparison with the previous work lie on the technical side. They
can be explained by considering the algebraic group whose real
points is $\PO(n,1)^\circ$. First, its algebraic simply
connected  covering is  a $4$-covering  when $n$  is odd  (while
it is only a  double covering  when  $n$ is even). Secondly, when
$n$ is odd, $\PO(n,1)^\circ$  is of type $\Dn$, for which there
exist outer  forms (for $n$  even the type  is $\Bn$ and all
forms of this type are inner).  In particular, the triality forms
$^{3,6}\Dn_4$ come out while dealing with $n = 7$. However, our
work shows that the orbifold $O^7_0$ is not defined by a triality
form. In order to deal with this issues we need to analyze in
detail the structure of the algebraic groups of type $\Dn$
defined over number fields. More precisely, compact hyperbolic
arithmetic orbifolds can be constructed using any totally real
number field different from $\Q$.

\medskip

Geometric  methods  can  be used  in  small dimensions   ($n   \le  9$)   to   determine  non-compact   hyperbolic $n$-orbifolds   of  the smallest  volume   (including non-arithmetic orbifolds) \cite{Mey, H, KH}. All these minimal orbifolds are known to be arithmetic.   The following  theorem  is consistent  with the  results of~\cite{H} and gives new results in higher dimensions:

\begin{Theorem}
  For each  dimension $n = 2r-1  \ge 5$, there is  a unique orientable
  non-compact arithmetic hyperbolic $n$-orbifold $O^n_1$ of the smallest
  volume. Its volume is given by:
  \begin{enumerate}[(1)]
  \item if $r \equiv 1 \;{\mathrm (mod\ 4)}$:
    \begin{eqnarray*}
      \vol(O^n_1) &=&  \frac{1}{2^{r-2}} \; \zeta(r) \; \prod_{i=1}^{r-1} \frac{(2i -1)!}{(2
        \pi)^{2i}} \zeta(2i);
    \end{eqnarray*}

  \item if $r \equiv 3 \;{\mathrm (mod\ 4)}$:
    \begin{eqnarray*}
      \vol(O^n_1)  &=&  \frac{(2^r  -1)  (2^{r-1}
        -1)}{3   \cdot  2^{r-1}}   \;  \zeta(r)   \;  \prod_{i=1}^{r-1}
      \frac{(2i -1)!}{(2 \pi)^{2i}} \zeta(2i);
    \end{eqnarray*}

  \item if $r$ is even:
    \begin{eqnarray*}
      \vol(O^n_1) &=&  \frac{3^{r-1/2}}{2^{r-1}} \;
      L_{\ell_1|\Q}\!(r) \; \prod_{i=1}^{r-1}
      \frac{(2i -1)!}{(2 \pi)^{2i}} \zeta(2i),
    \end{eqnarray*}
    where  $\ell_1 = \Q[\sqrt{-3}]$.
  \end{enumerate}

\end{Theorem}

The remarks made above about the proof of Theorem~1 remain valid
for Theorem~2, except some specific issues. Thus, by Godement's
compactness criterion the field of definition of a non-compact
arithmetic hyperbolic $n$-orbifold is known to be
the field of rational numbers $\Q$.
Also, the issue of triality forms does not appear here: non-compact
hyperbolic orbifolds are never related
to the triality (see Section~\ref{data-from-admissibility}).
Regardless of this few simplifications the proof of Theorem~2
requires essentially all the machinery which is used for
Theorem~1.

The results of Hild in \cite{H} which we mentioned above allow us
to compare the formulas from Theorem~2 with what was obtained by
geometric methods for small dimensions. We note that the results
of \cite{H} are not limited to orientable orbifolds. In particular, it is
proved there that the smallest volume $n$-orbifold for $n=5,7,9$
is non-orientable and unique. This shows that the orientable
double covers of
these orbifolds correspond to our $O_1^n$ ($n=5,7,9$).
Their volumes coincide with the formulas given in Theorem~2.

By evaluating the volumes in Theorems~1 and 2 and comparing their
asymptotic growth for large $n$ we obtain the following theorem of an
independent interest. It was already known for even dimensions
(see \cite[Sec.~4.5]{B1}), and our work allows us to extend it for all
sufficiently large $n$. This result gives further support to a
general conjecture stated in \cite[Add.~1.6]{B1}.

\begin{Theorem}
For every dimension $n\ge 5$, the minimal volume of a non-compact arithmetic hyperbolic $n$-orbifold
is smaller than the volume of any compact arithmetic hyperbolic $n$-orbifold. Moreover, the ratio
between the minimal volumes $\vol(O^n_0)/\vol(O^n_1)$ grows super-exponentially with $n$.
\end{Theorem}

It is interesting to point out a counterintuitive feature of this result: certain orbifolds
which are non-compact and thus have infinite cusps appear to be smaller than any compact
arithmetic orbifold of the same dimension. We can recall that for $2\le n \le 4$ the situation
is different, here the minimal volume of a non-compact hyperbolic $n$-orbifold is larger
than the volume of the smallest compact arithmetic $n$-orbifold (see \cite{CF}, \cite{Mey},
\cite{B1} and \cite{KH}). The same result holds true for $n=2$ and $3$ without arithmeticity
assumption. We conjecture that in higher dimensions one can also remove the assumption of
arithmeticity and that our Theorem~3 should apply to all hyperbolic $n$-orbifolds.

\medskip

We now give a brief outline of the paper. Section~\ref{section:prelim} provides some preliminary results which are used along the lines. In Section~\ref{section:candidates} we construct principal arithmetic subgroups $\Lambda_0$ and $\Lambda_1$ whose normalizers appear to be the natural candidates for the arithmetic subgroups of minimal covolume. This construction is motivated by Prasad's volume formula introduced in Section~\ref{sec:vol formula}. In Section~\ref{section:index-estimation} we give a bound for the index of a principal arithmetic subgroup in its normalizer. This bound is essential for our argument, its derivation is based on a detailed study of Galois cohomology. The results of Section~\ref{section:index-estimation} are used in Section~\ref{sec:index-candidates} to compute the precise indices for our groups $\Lambda_0$ and $\Lambda_1$. Uniqueness of the groups $\Lambda_0$, $\Lambda_1$ and their normalizers is established in Section~\ref{section:uniqueness}.
The proof is based on a reformulation of the problem in terms of Galois cohomology which allows us to apply the methods of the class field theory together with the results of the previous sections to show the uniqueness.
After all these preparations we are ready to prove the main theorems. This is done in Sections~7--10. The idea is first to use the general estimates and reduce the list of the possible candidates for the minimal volume. Then we improve the estimates using some specific properties of the candidates from the list and reiterate the procedure. After several steps of this optimization process we reduce the number of possibilities to a few candidates. A careful analysis of these remaining cases shows that the minimum is always attained on the groups which were constructed in Section~\ref{section:candidates}. This finishes the proof of Theorems~1 and 2. The final Section~\ref{section:growth} is dedicated to the study of growth of the minimal volume and the proof of Theorem~3. The results obtained here have some interesting applications. In particular, the method of \cite{ABSW} allows to apply these results together with \cite{B1} to the problem of classification of arithmetic hyperbolic reflection groups.

\medskip

\noindent
{\em Acknowledgements.} We would like to thank Ruth Kellerhals,
Gopal Prasad,  Alireza Salehi Golsefidy and Thomas Weigel for
their interest in this work and helpful discussions. We would
also like to thank the referee for their comments and
suggestions.


\section{Preliminaries on volumes and arithmetic subgroups}
\label{section:prelim}
\subsection{}\label{hyperbolic}

The group of isometries of the hyperbolic $n$-space
$\isom$ is isomorphic to the real Lie group $\PO(n,1)$.  The
group of orientation preserving isometries $\isomp$ is
isomorphic to  the identity  component $\bH =  \PO(n,1)^\circ$, which
can be further identified  with the matrix group $\SO(n,1)^\circ$. The
simply  connected  covering  of  $\bH$  is the  spinor  group  $\gH =
\Spin(n,1)$.  Any   orientable  hyperbolic  $n$-orbifold  is  a
quotient $\Hy^n/\overline{\Gamma}$ for some lattice $\overline{\Gamma}$ in $\bH$. We can
consider the full  inverse image $\Gamma$ of $\overline{\Gamma}$  in $\gH$.
A lattice in $\gH$ can be identified with such
an inverse image whenever it contains the center of $\gH$.

Let $\LieH$ be the Lie algebra of $\gH$ (and $\bH$). Each Haar
measure on $\gH$  or $\bH$ corresponds to a  multilinear form of the highest
degree on $\LieH$. If $K$ denotes a maximal compact subgroup of $\bH$
with Lie  algebra $\mathfrak{k}$, then $\Hy^n$ can  be identified with
$\bH/K$. Furthermore, the tangent space of $\Hy^n = \bH/K$ at the point
$K$ can be identified with a specially chosen subspace $\mathfrak{p} \subset \LieH$
such that   $\LieH   =    \mathfrak{k}   \oplus   \mathfrak{p}$.    It   is
known~\cite[Ch. V, \S 2]{Helg} that $\mathfrak{k}
\oplus i \mathfrak{p}$ is the (real) Lie algebra of the groups $\Spin(n+1)$
and $\SO(n+1)$, and this gives a duality between   $\Hy^n$ and the
$n$-sphere of curvature  one. Thus for any Haar  measure $\nu$ on $\gH$
there is a corresponding   Haar   measure  on   $\Spin(n+1)$, which we
will denote by $\nu^*$. Moreover, we denote by the
same symbol two Haar measures on $\gH$ and $\bH$ that correspond to a same
multilinear form on $\LieH$. We adopt the same convention for the measures on 
$\Spin(n,1)$ and $\SO(n+1)$.

Let  $\omega^1$   be  the   multilinear   form  on  $\mathfrak{p}$
corresponding to the left-invariant  measure $\vol$ and let $\omega^0$ be the
form on $\mathfrak{k}$ normalized by $\int_K \omega^0 = 1$.  Then the Haar measure
$\mu_\Hy$ given by  $\omega^0 \wedge \omega^1$ is such  that for any
lattice $\overline{\Gamma}$ in $\bH$ we have:
\begin{eqnarray*}
	\vol(\Hy^n /\overline{\Gamma}) & = & \mu_\Hy ( \bH /\overline{\Gamma}).
\end{eqnarray*}
Note that if $\Gamma$ contains the center of $H$, we
have $\mu_\Hy (\gH/\Gamma) = \mu_\Hy(\bH/\overline{\Gamma})$.
By the duality we can also check that $\mu^*_\Hy(\SO(n+1))$ is equal to the usual
volume of the unit $n$-sphere.

In this article $\mu$ stands for the Haar measure on $\gH$ normalized
by $\mu^*(\Spin(n+1)) = 1$.  This
is  the  normalization of  the  measure  used by Prasad in~\cite{P}. Since the
covering $\Spin(n+1) \to \SO(n+1)$ is of degree two, we have $\mu^*(\SO(n+1)) =
1/2$. Hence using the formula for the volume of the unit $n$-sphere in $\R^{n+1}$,
in the case of $n$ odd, we get:
\begin{eqnarray}
  \mu_\Hy              & = &
  \frac{4 \pi^\frac{n+1}{2}}{\left(\frac{n-1}{2}\right)!} \; \mu .
\end{eqnarray}
This allows us to work with $\mu$ instead of $\vol$.

\subsection{}\label{arithmetic-sbgp-def}

Let $\G$ be an absolute simple, simply connected algebraic group  defined over a number
field $k$.  We denote  by $V_f$ (resp.  $V_\infty$) the set  of finite
(resp. infinite) places of $k$.
A collection $\cP =  (\Pa_v)_{v\in V_f}$ of parahoric subgroups
$\Pa_v   \subset  \G(k_v)$ is called \emph{coherent} if
the product $\prod_{v \in V_f} \Pa_v$
is an open subgroup of the finite ad\`ele group $\G(\A_f(k))$
(we refer to e.g. \cite[Sec.~0.5]{P} for the definition of parahoric
subgroups). In this case
the intersection $\Lambda_\cP = \G(k) \cap \prod_{v \in V_f} \Pa_v$ is called
a \emph{principal arithmetic subgroup} of $\G(k)$.
Any subgroup  of $\G(k \otimes_\Q \R)$ that is commensurable with some principal arithmetic
subgroup is called an {\em arithmetic subgroup}.

By the Bruhat-Tits theory, the conjugacy class
(referred to as the \emph{type}) of a parahoric subgroup
$\Pa_v$ in $\G(k_v)$  corresponds  to a  subset $\theta_v$ of a basis $\Delta_v$
of the affine root system of $\G(k_v)$.
This way a principal arithmetic subgroup $\Lambda$ determines
a \emph{global type}
$$\theta  = \prod_{v  \in  V_f} \theta_v  \subset  \prod_{v \in  V_f}
\Delta_v . $$

An algebraic group $\G$ as above will be called {\em admissible} if there 
exists a continuous surjective homomorphism $\phi: \G(k\otimes_\Q\R) \to  \gH$ whose
kernel  is compact.  This implies  that $\G$  is simply  connected and
there  exists $v_0\in V_\infty$  such that  $\G(k_{v_0})\simeq\gH$ and
$\G(k_{v})$ is compact for  all $v\in V_\infty \setminus \{v_0\}$.
In particular, $k$ has to be a totally real field. Given an arithmetic
subgroup $\Gamma \subset \G(k \otimes_\Q \R)$ for some admissible group $\G$ defined over $k$,
its image $\phi(\Gamma)$ is a lattice in $\gH$. All
subgroups of $\gH$ which are commensurable with such $\phi(\Gamma)$ are called {\em arithmetic subgroups of $\gH$
defined over $k$}. Arithmetic subgroups of $\gH$ of the form $\phi(\Lambda)$ for some principal arithmetic
subgroup $\Lambda$ of $\G$ will be called the {\em principal arithmetic subgroups of $\gH$}.

\medskip

From now on $\G$  will always be an admissible group
defined over $k$ and of absolute rank $r \ge 3$,
$v_0 \in V_\infty$  denotes the archi\-medean place at which $\G$
is isotropic, and $\Delta_v\; (v \in V_f)$ denotes a chosen basis of the affine
root system of $\G(k_v)$.

\subsection{}
\label{data-from-admissibility}

The group $\G$ is of absolute type $\Dn_r$ (where $n = 2r -1$) and
we will often write  $\Dn_{2m}$ or $\Dn_{2m+1}$ to distinguish
between $r$ even or odd.
It follows from the classification of
algebraic groups  \cite{T1} and Godement's  compactness criterion that
$\G$ gives rise to
cocompact lattices in $H$ if and only if $k \neq \Q$ (see \cite[Secs. 1-2]{LM}
 for details).

By the admissibility condition, $\G(k \otimes_\Q \R)$ contains exactly one factor
$\Spin(n,1)$ and the rest of the factors, if any, are isomorphic to
$\Spin(n+1)$. By the classification of semisimple algebraic groups
(cf.~\cite[\S 3]{T1}) the group $\Spin(n,1)$ (resp. $\Spin(n+1)$) is an inner
form exactly when the discriminant of the standard quadratic form of signature $(n,1)$
(resp. the standard positive definite form) has the same sign as the
discriminant of the split form in $n+1$ variables, which equals
$(-1)^{r}$.
For even $r = 2m$ we conclude that $\Spin(n,1)$ is outer. This forces $\G$
to be an outer form (of possible types $^2\Dn_{2m}$ or the
exceptional \emph{triality forms} $^{3,6}\Dn_4$). For odd $r = 2m+1$
and $k \neq \Q$ the group $\G(k \otimes_\Q \R)$ must contain at least
one compact factor $\Spin(n+1)$, which is then an outer form. In this case again $\G$ is an outer
form (of the only possible type $^2\Dn_{2m+1}$). For the remaining case
$k =\Q$ with odd rank, both types $^1\Dn_{2m+1}$ and $^2\Dn_{2m+1}$ are possible.

\subsection{}
\label{sec:field-extension-l}

A  certain  field  extension  $\ell|k$  attached to  $\G$  will  be  used
throughout the paper.
For $\G$ of type different from $^6\Dn_4$,
the field $\ell$ is defined as the minimal extension
of $k$ such that $\G$ becomes an
inner form over $\ell$. In particular, except for the type $^3\Dn_4$, we have
$[\ell:k]\le 2$. For $\G$ of type $^6\Dn_4$ we let
$\ell$ to be a fixed extension of $k$ of degree $3$ which is contained in the
minimal extension of $k$ of degree $6$ over which $\G$ is inner. For
the types $^{3,6}\Dn_4$  we have $[\ell:k] = 3$.  By abuse of language
the field $\ell$  will be called in this  paper the \emph{splitting field}
of $\G$.

Let $(s_1,s_2)$ be the \emph{signature} of the field
$\ell$, i.e. $s_1$ (resp. $s_2$) is the number
of real (resp. complex) places of $\ell$. In
the following proposition we collect some useful information about
this field.

\begin{prop}
  \label{prop:facts-about-l}
  Let $\G|k$ be an admissible group of the absolute rank $r$ and let $d = [k:\Q]$.
  We suppose that $\G$ is an outer form.
  \begin{enumerate}
  \item If $r = 2m+1$, then the signature of $\ell$ is given
    by $(s_1,s_2) = (2,d-1)$.
  \item If $r = 2m$ and $\G$ is not
    of type $^{3,6}\Dn_4$ , then $(s_1,s_2) = (2d-2,1)$.
  \item When $k \neq \Q$ the field $\ell$ has at least one real
    place and the only roots of unity in $\ell$ are
    $\pm 1$.
  \end{enumerate}
\end{prop}
\begin{proof}
  The field $\ell$, except for the case $^6\Dn_4$,
  is the minimal splitting field of the quasi-split
  inner form of $\G$. This implies that for $v \in V_\infty$,
  if $\G$ is an inner form over $k_v$, the field $\ell$
  is a subfield of $k_v$ and hence for any
  place $w$ of $\ell$ extending $v$ the field $\ell_w$ is real. If $\G$ is an outer
  form over $k_v$ the opposite holds: there must exist $w|v$
  for which $\ell_w$ is complex. The two first assertions now follow from
  the description of the type of $\Spin(n,1)$
  and $\Spin(n+1)$ given in Section~\ref{data-from-admissibility}.

  If $k\neq\Q$ and the type is different from $^{3,6}\Dn_4$,
  the previous consideration proves also that $\ell$ must
  have at least one real place. For the cases
  $^{3,6}\Dn_4$ we have $[\ell:k]=3$ and then for each
  $v \in V_\infty$ we must have at least one real place $w|v$ of $\ell$,
  since $[\ell:k] = \sum_{w|v}[\ell_w:k_v]$.
  Thus in each of the cases $\ell$ has at least one real place, and hence the
  only roots of unity in $\ell$ are $\pm 1$.
  This proves the third assertion.
\end{proof}

\subsection{}\label{sec:vol formula}

The covolume of a principal arithmetic subgroup $\Lambda \subset
\G(k_{v_0}) \simeq \gH$ corresponding  to a coherent collection $\cP
= (\Pa_v)_{v \in V_f}$ can be computed
using Prasad's volume formula~\cite[Theorem~3.7]{P}. In our case
Prasad's formula gives
\begin{equation}
\mu(\gH/\Lambda) = \D_k^{\frac{2r^2-r}2} \left(\frac{\D_\ell}{\D_k^{[\ell:k]}}\right) ^{\frac12(2r-1)}\;
C(r)^{[k:\Q]} \; \cE(\cP),
\end{equation}
where $\ell$ is the field introduced in Section~\ref{sec:field-extension-l}, $\D_K$
denotes the absolute value of the discriminant of a number field
$K$ (in this paper we will call $\D_K$ briefly the \emph{discriminant}
of $K$), $\cE(\cP) = \prod_{v\in  V_f} e_v(\Pa_v)$ is an Euler product
of local factors which are determined by the structure of $(\Pa_v)_{v\in V_f}$, and
the constant $C(r)$ is defined by
$$
C(r) = \frac{(r-1)!}{(2\pi)^{r}}\;\; \prod_{i=1}^{r-1} \frac{(2i-1)!}{(2\pi)^{2i}} .
$$

\subsection{}\label{sec:local-factor}

The local factors $e(\Pa_v)$ can be computed using the Bruhat-Tits theory. An
extensive account of its main concepts and results is given in \cite{T2}.

Let  $k_v$  be a  nonarchimedean  local  field,  $\kun_v$ its  maximal
unramified extension, and $\Pa_v$ be a parahoric
subgroup  of   $\G(k_v)$.  Following  \cite[Sec.~3.7]{P},
$e(\Pa_v)$ is given by the formula
\begin{equation*}
e(\Pa_v) = \frac{q_v^{(\dm\Mv + \dm\cMv)/2}}{\#\Mv(\fv)},
\end{equation*}
where $q_v$ is the degree of the residue field $\fv$ of $k_v$,
and $\Mv$, $\cMv$ are certain connected reductive
$\fv$-groups associated to $\G$ and its quasi-split inner form
which are defined in \cite[Sec.~2.2]{P}.
There exists a minimal finite (possibly empty) subset $T\subset V_f$,
such that for all $v\in  V_f\setminus T$, the  group $\G$ is quasi-split
over $k_v$ and $\Pa_v$ is a hyperspecial (resp. special) parahoric
subgroup if $\G$ splits over $\kun_v$ (resp.  does not
split over $\kun_v$).
For $v\in  V_f\setminus T$  the group  $\Mv$ is  isomorphic to
$\cMv$ over $\fv$.

We assume till the end of Section~\ref{sec:local-factor} that $\G$ is not of
type  $^{3,6}\Dn_4$.   If  $\G$  splits  over  $\kun_v$   and  $v  \in
V_f\setminus T$, then $\Mv \simeq \cMv$ is a connected semisimple
$\fv$-group of type $^1\Dn_r$ or $^2\Dn_r$ according
with $\G$ being split or non-split  over $k_v$. If $\G$ does not split
over $\kun_v$ and $v \in V_f \setminus T$, then $\Mv \simeq \cMv$ is of type $\Bn_{r-1}$.
Hence for $v\in V_f\setminus T$, we have (see \cite[Table~1]{Ono}):
\begin{gather*}
\dm\Mv =  \dm\cMv = 2r^2 - r,\\
\#\Mv(\fv) = q^{r(r-1)}(q^r\mp 1)\prod_{i = 1}^{r-1}(q^{2i}-1),
\end{gather*}
if  $\G$ splits  over $\kun_v$,  the sign  $\mp$ being
negative (resp.  positive) if $\G|k_v$ is inner (resp. outer); and
\begin{gather*}
\dm\Mv =  \dm\cMv = 2(r-1)^2 + r-1, \nonumber\\
\#\Mv(\fv) = q^{(r-1)^2}\prod_{i = 1}^{r-1}(q^{2i}-1), \label{eq Mv}
\end{gather*}
if $\G$ does not split over $\kun_v$.

We can bring these cases together using the field $\ell$ defined
in Section~\ref{sec:field-extension-l}.
Namely, assuming for a moment that $T = \emptyset$, we obtain
\begin{equation*}
\prod_{v\in V_f} e(\Pa_v) = \zeta_k(2)\zeta_k(4)...\zeta_k(2r-2) L_{\ell|k}(r).
\end{equation*}
Now, for each $v\in T$ we define a so-called \emph{lambda factor}:
\begin{eqnarray}
 \lambda_v &=& e(\Pa_v) \; \cdot \; \frac{\# \cMv(\fv)}{q_v^{\dim(\cMv)}}.
\end{eqnarray}
 Then
\begin{equation}
\cE(\cP) = \zeta_k(2)\zeta_k(4)...\zeta_k(2r-2) L_{\ell|k}(r) \prod_{v\in T}
\lambda_v,
\end{equation}
moreover, by \cite[Prop.~2.10]{P},
\begin{align}
  \label{eq:bound-for-lambda-factor}
\lambda_v &\ge
(q_v+1)^{-1}q_v^{r_v+1} (1 \mp q_v^{-r_v})\prod_{i = 1}^{r_v-1}(1 - q_v^{-2i}) \nonumber\\
&\ge \frac{2}{3} \left(\frac{3}{4} q_v \right)^{r_v},
\end{align}
where $r_v$ is the $\kun_v$-rank of $\G$.

\subsection{}
\label{sec:minimal-volume-implies-max}
It is well known that every maximal arithmetic
subgroup $\Gamma$ of $\gH = \G(k_{v_0})$ is a normalizer
$N_{\gH}(\Lambda)$ of some
principal arithmetic subgroup $\Lambda$ (see e.g. \cite[Prop.~1.4]{BP}).
In~\cite{R} (see also~\cite{CR}) Rohlfs gave a characterization of
principal arithmetic subgroups whose normalizers are maximal by defining
a notion of ``$\mathcal{O}$-maximality''.
Arithmetic subgroups  of minimal covolume are  necessarily maximal. In
order to determine these subgroups  we will first look at the
principal arithmetic subgroups which respect Rohlfs' criterion and
have small covolume (Section~\ref{section:candidates}), and
later on will consider their normalizers.


\section{\texorpdfstring{Construction of $\Lambda_0$ and $\Lambda_1$}{Construction}}
\label{section:candidates}
\newcommand{\disc}{\delta}

\subsection{}
\label{sec:examples_intro}

In this section we give a construction of the principal arithmetic
subgroups $\Lambda_0$ and $\Lambda_1$, whose normalizers
define the minimal orbifolds from Theorems 1 and 2.
These will be subgroups of admissible spinor groups
$\Spin_f$ for some quadratic forms $f$.
The group $\G=\Spin_f$ needs to fulfill the admissibility condition
of Section~\ref{arithmetic-sbgp-def}. This is precisely the case when $f$ is
defined over a totally real number field $k$, has signature
$(n,1)$ at an infinite place $v_0$ of $k$ and is
anisotropic for all $v \in$ \mbox{$V_\infty \setminus \left\{v_0\right\}$}.
We call such $f$ an \emph{admissible quadratic form}.
Arithmetic subgroups of such groups $\G$ are sometimes referred to as arithmetic
lattices of the \emph{first type}. Such a lattice is cocompact if and only if
$k  \neq   \Q$  (see Section~\ref{data-from-admissibility}).  All non-cocompact
arithmetic subgroups  are of the first type \cite[Secs.  1-2]{LM}, i.e.
every  admissible  $\Q$-group  is  of  the  form  $\Spin_f$  for  some
admissible quadratic form $f$ defined over $\Q$.

\subsection{}
\label{sec:examples_alg-struct}

Let $(\V,f)$  be a $2r$-dimensional  quadratic space over a  field $k$
and let $\h$ be the two-dimensional isotropic quadratic space.
Then $\G=\Spin_f$ is quasi-split over $k$ if and only if
$(\V,f) \simeq \h^{r-1} \oplus D$
for a  two-dimensional space $D$, and $\G$ is split or not
according to $D \simeq \B{h}$ or not (see \cite[\S 23.4]{B-book}
where the result is explained for $\SO_f$).

We remind that when $k$ is a number field,
quadratic forms of a given degree are
classified by looking at every local extension $k_v$ of $k$
(\emph{the Hasse-Minkowski theorem}).
Quadratic forms over non-archimedean extensions are
classified by the discriminant (denoted by $\disc$) and the \emph{Hasse symbol}.
We refer to O'Meara's monograph~\cite{OM} for the facts about quadratic forms.
By the Hasse symbol we mean the normalization of the Hasse invariant which is
used in this book.

The group $\G = \Spin_f$ over a non-archimedean local field $k_v$ has the following structure:
\begin{enumerate}[(1)]
\item $\G$ is a split inner form if and only if $\disc (f)$
  is equal to $\disc (\B{h}^r) = (-1)^r$ up to an
  element of $(k_v^*)^2$.
  In this case $\G$ is
  split exactly when $f$ has the same Hasse symbol
  as $\h^r$.

\item For other values of $\disc(f)$, $\G$ is quasi-split
  but non-split and $\disc(f) \in k_v^*/(k_v^*)^2$ determines uniquely the $k_v$-isomorphism class
  of $\G$.
\end{enumerate}

\subsection{}
\label{sec:on-skew-hermitian}

Recall  that in  odd  hyperbolic dimension,  besides the  construction
based on quadratic forms, there exist
(cocompact) arithmetic lattices defined by skew-hermitian forms over division quaternion
algebras (see \cite[Sec.~2]{LM} for the details of this construction).
Moreover, Cayley algebras produce examples
of arithmetic subgroups related to triality forms $^{3,6}\Dn_4$.
If  we exclude  this latter  case (which  concerns only the hyperbolic dimension
seven), we  can use the language  of algebras with  involution to unify
the  description of  admissible  groups (we  refer to~\cite{KMRT} for the
notions related to this subject).  Namely, all non-trialitarian $k$-groups which
are admissible  in the sense of Section~\ref{arithmetic-sbgp-def}  have (up to
$k$-isomorphism) the  form $\Spin(A,\sigma)$,  where $A$ is  a central
simple  algebra  and  $\sigma$  is  an  orthogonal  involution.   More
precisely, $A$ is of the form $\End_R(M)$ for a quaternion algebra $R$ over $k$ and
an $R$-module $M$.
If the admissible group is defined by a skew-hermitian form $h$, then $R$
corresponds to the division algebra over which $h$ is defined
\cite[Theorem 4.2]{KMRT}. On the other hand,
the algebra $R$ \emph{splits}, i.e.  $R
\simeq \End_k(k^2)$, exactly when $A$  has the form $\End_k(\V)$ for a
quadratic space $(\V,f)$ with $\sigma$ being the associated adjoint
involution.  In this  case $\Spin(A,\sigma)$ naturally identifies with
the $k$-group $\Spin_f$. The following result can be considered as a
variation on the Hasse-Minkowski theorem.

\begin{lemma}
  \label{lem:uniqueness-G0}
  Let  $\G$  be  an   admissible  $k$-group  of  type  different  from
  $^{3,6}\Dn_4$.  Suppose that  for  each  place $v  \in  V_\infty\cup V_f$, 
  $\G$ is $k_v$-isomorphic to some group $\Spin$ of a quadratic form over
  $k_v$. Then  $\G$ is $k$-isomorphic to $\Spin_f$  for some quadratic
  form  $f$ over $k$.
\end{lemma}

\begin{proof}
  By Section~\ref{sec:on-skew-hermitian},   the   group    $\G$ has the
  form $\Spin(\End_R(M),\sigma)$ for some  quaternion algebra $R$ over $k$.
  By the hypothesis $R$  splits  over each  completion  $k_v$ of the field $k$.
  This implies   that   $R$   splits   already   over   $k$   \cite[Theorem~2.7.2]{McR}.
  So $\G$  has the  form $\Spin_f$  for  some quadratic form $f$.
  \end{proof}

\subsection{}
\label{sec:constr-G-0}

Let us now construct an admissible group $\G_0$ which will give rise to a
cocompact arithmetic subgroup of small covolume.
In this case we have to consider real number fields $k$
different from $\Q$. To get small covolume,
by Section~\ref{sec:vol  formula}, we  need to  consider fields  $k$  and $\ell$
whose discriminants have small  absolute values. The smallest possible
value for $\D_k$ is attained
for $k = k_0 =\Q[\sqrt{5}]$. The smallest discriminant of a field
$\ell$ which satisfies the conditions of Proposition~\ref{prop:facts-about-l}
is $\D_\ell =  275$. It corresponds to a quartic field with a
defining polynomial $x^4 - x^3 + 2x -1$. This pair $(k_0,\ell_0)$ comes  out, for example,
with the following quadratic form:
\begin{equation*}
  \label{eq:form-cocompact}
  f_0(x_0,\dots,x_n) = \left((-1)^{r} \cdot3   - 2 \sqrt{5} \right) x_0^2 \;
+\; x_1^2 \;+\; \cdots \;+\; x_n^2.
\end{equation*}

\begin{prop}
  \label{prop:G-compact-quasi-split}
  The group $\G_0 = \Spin_{f_0}$ is quasi-split at every finite place
  $v \in V_f$. It is the unique admissible group attached to
  $(k_0,\ell_0)$ which respects this property.
\end{prop}

\begin{proof}
	The ideal of $\mathcal{O}_{k_0}$ generated by $(-1)^r \disc(f_0)$
  is prime. It corresponds to the unique finite place
  $v_\mathrm{r}$ of $k_0$
  at which the extension $\ell_0|k_0$ ramifies. At this ramified place $\G_0$
  cannot be inner, and from Section~\ref{sec:examples_alg-struct}
  it follows that it has to be quasi-split.
  Let $v_2$  be the finite place  of $k_0$ corresponding  to the prime
  ideal generated  by $2$.  Then for every  finite place  $v$ different
  from $v_2$ and $v_\mathrm{r}$ the fact that the discriminant $\disc(f_0)$ is a
  unit at $v$ implies that the Hasse invariant is trivial and thus $\G_0$
  is quasi-split over $k_v$.
  For the remaining place $v_2$ one can
  check that the equations $x^2 = \pm \disc(f_0)$ have no
  solution modulo~$(8)$ in the ring of integers of $k_0$, proving that
  $\pm \disc(f_0)$ cannot be squares in the local field $(k_0)_{v_2}$.
  Hence by Section~\ref{sec:examples_alg-struct},
  $\G_0$ is quasi-split over $(k_0)_{v_2}$ as well.

  To prove  the uniqueness result  we first note that  a quasi-split
  group of type $^1\Dn_r$ or $^2\Dn_r$ is necessarily given by a quadratic form
  (see \cite[Example~27.10]{KMRT}). Since admissibility implies that the group
  is also defined by quadratic forms over the archimedean completions (namely, by
  the  forms  of  signature  $(n,1)$  and  $(n+1,0)$),  we  can  apply
  Lemma~\ref{lem:uniqueness-G0}  to deduce  that  an admissible  group
  respecting the property stated for the group $\G_0$ in the proposition must be
  of  the type  $\Spin_f$, for  some  quadratic form  $f$ defined  over
  $k_0$.  Specifying that $\Spin_f$ has the same splitting field $\ell_0$ as the
  group $\G_0$ implies that $\Spin_f \simeq \G_0$ over each completion
  of $k_0$.  By the Hasse-Minkowski theorem we conclude that such a group
  $\Spin_f$ is $k_0$-isomorphic to $\G_0$.
\end{proof}

\subsection{}
\label{sec:constr-Gamma-0}
In order to construct a principal arithmetic
subgroup of small covolume we need to consider
parahoric subgroups of the maximal volume.
When $\G$ is quasi-split over $k_v$  and splits over its maximal
unramified extension there
exist hyperspecial parahoric subgroups in $\G(k_v)$
and those are of the maximal volume.
If $\G$ does not respect these hypothesis, the
parahoric subgroups of maximal volume
are the special ones (see \cite[Sec.~3.8.2]{T2}).

Proposition~\ref{prop:G-compact-quasi-split} implies the existence of a coherent
collection $\cP$ of parahoric subgroups $\Pa_v \subset \G_0(k_v)$
which are hyperspecial for every place $v$ of $k = k_0$ unramified in $\ell_0$
and special (but not hyperspecial) for the unique ramified
place $v_\mathrm{r}$ corresponding to the ideal $(3-(-1)^r 2 \sqrt{5})$.
Let us denote by $\Lambda_0$
the principal arithmetic subgroup attached to $\cP$.
Then by Sections~\ref{sec:vol formula} and \ref{sec:local-factor}
its covolume is given by:
\begin{equation}
  \label{eq:vol_compact}
  \mu(H/\Lambda_0) = 5^{r^2 -r/2} 11^{ r - 1/2}\;\;C(r)^2\;\;
\underbrace{ \zeta_k(2)\zeta_k(4) \cdots \zeta_k(2r-2) L_{\ell|k}(r)}_{(*)}.
\end{equation}

Similarly  to  the  way  it  was  done  in  \cite[proof  of  Thm.~4.1,
p. 760]{B1}, we can show
that the product~$(*)$ is bounded by 2.
Hence we get
\begin{equation}
  \label{eq:bound_compact}
  \mu(H/\Lambda_0) \quad <\quad  2\;\cdot \; 5^{r^2 -r/2}\;\; 11^{ r - 1/2}\;\;  C(r)^2\;\;.
\end{equation}
For small $r$  we will use a better estimate for the product
$(*)$:  We can check  using  Pari/GP (which  can  compute  such  an
expression  up to a given precision) that this  product  is bounded  by
$1.17$ when $r \le 16$. Thus for these $r$ we have a better bound:
\begin{equation}
  \label{eq:bound_compact-1.2}
  \mu(H/\Lambda_0) \quad <\quad 1.17\;\cdot\;5^{r^2 -r/2}\;\; 11^{ r - 1/2}\;\; C(r)^2\;\;.
\end{equation}
These bounds will be used in the proof of Theorem~1 in Sections~\ref{section:proof-compact-odd} and
\ref{section:proof-compact-even}.

\subsection{}
\label{sec:non-compact_ex}

Now let us consider the non-compact case. It corresponds to  $k = \Q$.
In contrary to the preceding discussion here we have to distinguish
some cases according to the dimension $n = 2 r - 1$.

First consider $r$ even. By Section~\ref{data-from-admissibility},
the  field $\ell$ is a non-trivial extension of $\Q$, and
by Proposition~\ref{prop:facts-about-l} we
see that $\ell$ must be a quadratic imaginary field.
Among these fields the smallest discriminant is attained for
$\ell_1 = \Q[\sqrt{-3}]$. Let us take
$$
f_1(x_0, \dots, x_n) = -3 x_0 ^2 + x_1^2 + \cdots + x_n^2.
$$
For odd $r$ spinor groups of inner type come out. Here we can take
$$
f_1(x_0, \dots, x_n) = - x_0 ^2 + x_1^2 + \cdots + x_n^2.
$$
With the quadratic form $f_1$ defined this way (depending on $r$)
and  $\G_1 = \Spin_{f_1}$, the ``non-compact'' analogue of
Proposition~\ref{prop:G-compact-quasi-split} is given by the following
proposition.
\begin{prop}\
  \label{prop:G-local-non-compact}
  \begin{enumerate}
  \item If $r$ is even, $\G_1$ is quasi-split at every finite
    place. It is the unique admissible $\Q$-group with the splitting field
    $\Q[\sqrt{-3}]$ which has this property.
  \item  For  $r =  2m  +  1$ odd,  $\G_1$  is the unique  admissible
    $\Q$-group of type $^1\Dn_r$ which respects the conditions:
    \begin{enumerate}[(a)]
    \item if $m$ is even, $\G_1$ is split
      over $\Q_p$ for every prime $p$;
    \item if $m$ is odd, $\G_1$ is split over $\Q_p$ for every
      prime $p \neq 2$ and  $\G_1|\Q_2$
      is a split inner form but non-split.
    \end{enumerate}
  \end{enumerate}
\end{prop}
\begin{proof}
  Since the admissible $\Q$-groups  are all of the form  $\Spin_f$, in all
  the cases  the  uniqueness statement  follows  from the Hasse-Minkowski
  theorem  in  a   similar way as in the   proof  of
  Proposition~\ref{prop:G-compact-quasi-split}.  It remains to
  prove the specific properties of $\G_1$.

  When $r$ is even the situation is entirely
  similar to that in Proposition~\ref{prop:G-compact-quasi-split}:
  the Hasse invariant ensures that $\G_1$ is quasi-split
  at primes $p \neq 2,3$, and as $\ell_1|\Q$ is
  ramified at $3$, $\G_1$ is also quasi-split over $\Q_3$.
  The equation $x^2 = -3$ has no solutions
  in $\Q_2$ therefore $\G_1$ is quasi-split over $\Q_2$ as well.

  For $r$ odd, $\G_1$ is inner at every place
  of $\Q$. The Hasse invariant is trivial
  at each prime $p \neq 2$. Computing
  the Hasse invariant of $f_1$ and $\h^{r}$ over
  $\Q_2$ leads to the different results
  for $m$ even or odd, which give the second part
  of the proposition.
\end{proof}
This proposition implies the existence of a coherent
collection $\Cal{P}$ of parahoric subgroups $\Pa_v$ which
are all hyperspecial except for at most one place $v = (p)$ of $\Q$,
where we choose $\Pa_v$ to be special. For $r$ even this
exception is attached to the prime $p=3$ which is ramified in $\ell_1$.
For $r = 2m+1$: if $m$ is odd, $\Pa_v$ is not hyperspecial for $v=(2)$,
whereas for $m$ even all parahoric subgroups can be chosen to be hyperspecial.

Let $\Lambda_1$ be the principal arithmetic subgroup
of $H$ determined by $\Cal{P}$. In Table 1 we list
the covolumes of $\Lambda_1$ for the different cases.
\setlength{\extrarowheight}{5 pt}
\begin{table}[h]
  \label{tab:vol_nc}
  \centering
  \begin{tabular}[d]{|l|c|c|c|}
    \hline
    $r$ & splitting field $\ell$ & $ \mu(H/\Lambda_1)$ \\[3pt]
    \hline
    $r$ even & $\Q[\sqrt{-3}]$ & $3^{r -1/2} \;C(r)
    \zeta(2)\zeta(4) \cdots \zeta(2r-2) L_{\ell|\Q}(r)$
    \\[3pt]
    \hline
    $r \equiv 1 \; \mathrm{mod}\; 4 $ & $\Q$  &  $C(r)
    \zeta(2)\zeta(4) \cdots \zeta(2r-2) \zeta(r) $\\[3pt]
    \hline
    $r \equiv 3 \; \mathrm{mod}\; 4 $ & $\Q$  &  $ \lambda_{(2)} C(r)
    \zeta(2)\zeta(4) \cdots \zeta(2r-2) \zeta(r) $\\[3pt]
    \hline
  \end{tabular}
  \vspace{0.5cm}
  \caption{Covolume of $\Lambda_1$}
\end{table}
Note that in the last
case the value of $\lambda_{(2)}$ (see Section~\ref{sec:local-factor}) is
non-trivial. It can be computed using the description of the
reductive group $\overline{\mathrm{M}}_{(2)}$  attached to the
special parahoric subgroup $\Pa_{(2)}$.
This description is obtained from the local index of $\G(\Q_2)$,
which according to Tits \cite[Sec.~4.3]{T2} is
$^2\Dn'_r$  (resp.  $^2\mathrm{A}_3$  for $r=3$).   More  generally, we
can compute that if $\G(\Q_p)$  is of  this local  type, then  the lambda
factor coming from a special parahoric subgroup is given by:
\begin{eqnarray}
  \label{eq:lambda-factor-2}
  \lambda_{(p)} &=& \frac{\left( p^r-1 \right) \left(p^{r-1}-1\right)}{p+1}.
\end{eqnarray}

\subsection{}
\label{sec:final-rmk-cand}

Using Prasad's volume formula it can be shown without too much
effort that $\Lambda_0$ (resp. $\Lambda_1$) is of minimal
covolume    in    the     class    of    \emph{principal}    cocompact
(resp. non-cocom\-pact) arithmetic subgroups
of $H$. However, these groups are not necessarily maximal
arithmetic subgroups of $H$. Their normalizers
$$\Gamma_0 := N_{H}(\Lambda_0) \quad \mbox{ and }\quad \Gamma_1 := N_{H}(\Lambda_1)$$
are maximal (see Section~2.8)
and they will turn out to be the hyperbolic arithmetic lattices of minimal
covolume. In the next
section we will explain how the computation of the index
$[N_H(\Lambda):\Lambda]$  (with a principal  arithmetic subgroup $\Lambda$) can  be
carried out.  General estimates obtained there are used in the proofs of Theorems~1
and~2. In
Section~\ref{sec:index-candidates}   we  deal   specifically   with  the
computation of $[\Gamma_0:\Lambda_0]$ and  $[\Gamma_1:\Lambda_1]$.


\section{\texorpdfstring{Estimating the index $[\Gamma:\Lambda]$}{Estimating the index}}
\label{section:index-estimation}
We use the notation of Section 2. In particular,
$\G$ is an admissible group defined over a totally real number field $k$
and $\Lambda$ denotes a principal arithmetic subgroup of $\gH = \G(k_{v_0})$
with a normalizer $\Gamma$.

\subsection{}\label{sec:action-H1-on-parahoric}

Let $\CG$ be the center of $\G$ and $\bG$ be the adjoint group.
We have an exact sequence of $k$-isogenies:
\begin{equation*}\label{eq22}
1 \to \CG \to \G \stackrel{\phi}\to \bG \to 1.
\end{equation*}
This   induces   an  exact   sequence   in   Galois  cohomology   (see
\cite[Sec.~2.2.3]{PR}):
\begin{equation}\label{eq23}
\G(K) \stackrel{\phi}\to \bG(K) \stackrel{\delta}\to \HH^1(K, {\CG}) \to \HH^1(K,\G),
\end{equation}
valid for every field extension $K|k$.

The group $\G$ is simply connected hence
for   each  $v   \in   V_f$   we  have   $\HH^1(k_v,\G)   =  1$   (see
\cite[Theorem~6.4]{PR}).
The group $\bG(k_v)$ acts by conjugation on
the set of parahoric subgroups and the action of $\phi(\G(k_v))$ is trivial, hence
this induces a homomorphism
$$\xi_v:  \HH^1(k_v,\CG) \to  \Aut(\Delta_v).$$
The image of $\xi_v$ is denoted by $\Xi_v$, and
we can define the map
$$
\xi: \HH^1(k,\CG) \to \bigoplus_{v \in V_f} \Xi_v
$$
as a product of all $\xi_v$ (restricted to $\HH^1(k,\CG)$).
This material is explained in more detail in~\cite[Sec. 2]{BP}.

\subsection{}
\label{sec:Rohlfs-exct-sequ}

Given $E \le \HH^1(k,\CG)$, we denote
by $E_\xi$ its subgroup which acts trivially on each $\Delta_v$ and by
$E_\theta$ the subgroup which stabilizes the  type $\theta
= \prod_{v \in V_f} \theta_v $ attached to $\Lambda$
as in Section~\ref{arithmetic-sbgp-def}. We have
\begin{equation}
  \label{eq:24}
  \# E_\theta \le \# E_\xi \prod_{v \in V_f} \# \Xi_{\theta_v},
\end{equation}
where $\Xi_{\theta_v}$ is the subgroup of $\Xi_v$ stabilizing $\theta_v$.

Now by \cite{R}, \cite[Prop.~2.9]{BP}, we have an exact sequence
\begin{equation}\label{eq25}
1 \to \CG(k_{v_0}) / \left( \CG(k)\cap\Lambda \right) \to \Gamma/\Lambda \to \delta(\bG(k))'_\theta \to 1,
\end{equation}
where
\begin{eqnarray*}
\delta(\bG(k))' &=& \delta(\bG(k)) \cap (\delta \circ \phi)(\G(k_{v_0})), \text{\ and}\\
\delta(\bG(k))'_\theta &=& \delta(\bG(k))'\cap\HH^1(k,\CG)_\theta.
\end{eqnarray*}

In  our  case $\Lambda$,  being  a  principal  arithmetic subgroup  of
$\G(k)$, contains the center $\CG(k)$.
Also, $\CG(k_{v_0}) = \left\{ \pm 1 \right \}$
and using the description of $\CG$ in Section~\ref{sec:center-description} we can check that
if $\G$ is not of type $^{3,6}\Dn_4$, then $\CG(k) = \CG(k_{v_0})$.
This implies that for these types
$\Gamma/\Lambda  \simeq  \delta(\bG(k))'_\theta$.  For  $\G$  of  type
$^{3,6}\Dn_4$ we have $\CG(k) = 1$ and
the quotient $\Gamma/\Lambda$ is a $\Z/2\Z$-exten\-sion
of $\delta(\bG(k))'_\theta$.

\subsection{}
\label{sec:Lambda_m}

In \cite[Sec. 3]{BP} a principal arithmetic subgroup
$\Lambda^\mmax$  commensurable with $\Lambda$  is chosen  in such  a way
that  the  \emph{generalized index}  of  $\Lambda$ in  $\Lambda^\mmax$
respects the following inequality:
\begin{eqnarray}
  \label{eq:gener-index-for-mmax}
[\Lambda^\mmax : \Lambda] & = &
\frac{[\Lambda^\mmax : \Lambda^\mmax \cap \Lambda]}{[\Lambda : \Lambda^\mmax \cap \Lambda]}
 \ge \prod_{v \in V_f} \# \Xi_{\theta_v}.
\end{eqnarray}
We let $\theta^\mmax = (\theta_v^\mmax)$ be the global type of the
principal subgroup $\Lambda^\mmax$. By the construction in \cite{BP},
the types $\theta^\mmax_v$ are all special which
implies that for every $v \in V_f$ we have $\Xi_{\theta_v^\mmax} = 1$.
It follows that $\delta(\bG(k))'_{\theta^\mmax} = \delta(\bG(k))'_\xi$.
As this value depends only on the group $\G$, we conclude that the
covolume of the  normalizer $\Gamma^\mmax = N_H(\Lambda^\mmax)$ depends
only     on~$\G$.

Now we have
\begin{eqnarray*}
  \mu(H/\Gamma)    & = & \frac{\mu(H/\Lambda)}{[\Gamma:\Lambda]}  =
  \frac{\mu(H/\Lambda^\mmax)[\Lambda^\mmax : \Lambda]}{[\Gamma:\Lambda]}.
\end{eqnarray*}
By~\eqref{eq:gener-index-for-mmax}, it gives
\begin{eqnarray}\label{eq:bnd4.3.1}
  \mu(H/\Gamma)    & \ge &
  \frac{\mu(H/\Lambda^\mmax)}{[\Gamma:\Lambda]}\prod_{v \in V_f} \# \Xi_{\theta_v}.
\end{eqnarray}
From the other hand, by~\eqref{eq:24} and \eqref{eq25},
\begin{eqnarray}\label{eq:bnd4.3.2}
    [\Gamma:\Lambda]  & \le & [\Gamma^\mmax:\Lambda^\mmax] \prod_{v \in V_f} \# \Xi_{\theta_v}.
\end{eqnarray}
Combining inequalities \eqref{eq:bnd4.3.1} and \eqref{eq:bnd4.3.2} we obtain
\begin{eqnarray*}
  \label{eq:Lambda-m}
  \mu(H/\Gamma)    & \ge& \frac{1}{[\Gamma^\mmax:\Lambda^\mmax]}
  \mu(H/\Lambda^\mmax) =  \mu(H/\Gamma^\mmax).
\end{eqnarray*}
This means that $\Gamma^\mmax$ is of the smallest covolume
among the arithmetic subgroups attached to~$\G$. Its covolume
depends on $\#\delta(\bG(k))'_\xi$ and we will now
focus on the computation of this order.

In the rest of the article we
will  write $\Lambda =  \Lambda^\mmax$  to
indicate that  $\Lambda$ is  commensurable with $\Lambda^\mmax$ and has
a similar local  structure, i.e., $\theta^\mmax$ is in the
$\bG(\A_f)$-orbit of $\theta$ (cf. Section~\ref{sec:uniqueness-basic}). Similarly, the notation  $\Gamma^\mmax = \Gamma$
will be used for their normalizers. The equality is an abuse of
notation.  However, it is  clear that in this case $\Gamma$ (resp. $\Lambda$) has
the same covolume as $\Gamma^\mmax$ (resp. $\Lambda^\mmax$).

\subsection{}
\label{sec:center-description}
We remind here the description of the center $\CG$ from \cite[Table on
p. 332]{PR}.
Let $\nc=4$, $\epsilon = 1$ if $\G$ is of absolute type $\Dn_{2m+1}$,
and $\nc=2$, $\epsilon = 2$, otherwise.
The inner form of type $^1\Dn_{r}$ has as its center the group
$\mu_\nc^\epsilon$ of order $\nc^\epsilon$.
For the type $^2\Dn_{2m}$ the center $\CG$ is isomorphic to
$\Res_{\ell|k} (\mu_2)$, where $\Res_{\ell|k}$ denotes the
\emph{restriction of scalars} relative to the field extension $\ell|k$
defined in Section~\ref{sec:field-extension-l}.
In the remaining cases ($^2\Dn_{2m+1}$ and $^{3,6}\Dn_{4}$) the center
is isomorphic to the kernel  $\Res^{(1)}_{\ell|k}(\mu_\nc)$ of the \emph{norm
map} $\No_{\ell|k}: \Res_{\ell|k} (\mu_\nc) \to \mu_\nc$. Here we have
 the following exact sequence:
\begin{equation}
  \label{eq:center-exact-sq}
  1 \to \Res^{(1)}_{\ell|k}(\mu_\nc) \to \Res_{\ell|k} (\mu_\nc)
  \stackrel{\No_{\ell|k}}\to \mu_\nc \to 1.
\end{equation}

\subsection{}
\label{sec:cohom-center}
From Section~\ref{sec:center-description} we can deduce a description
of $\HH^1(k,\CG)$ for the admissible groups.
For the type $^1\Dn_{2m+1}$, we have $\CG = \mu_4$ and $\HH^1(k,\CG) = k^*/(k^*)^4$.
If $\G$ is of type $^2\Dn_{2m}$, then
$\HH^1(k,\CG) = \HH^1(k,\Res_{\ell|k} (\mu_2)) = \ell^*/(\ell^*)^2$.
Dealing with the remaining cases requires considering the cohomological exact
sequence associated to~\eqref{eq:center-exact-sq}, which gives us the
following exact sequence:
\begin{equation}
  \label{eq:cohom-center-ex-sq}
  1 \to \mu_\nc(k) / \No_{\ell|k}(\mu_\nc(\ell)) \to \HH^1(k,\CG) \to \ker\left(
    \ell^*/(\ell^*)^\nc \stackrel{\No_{\ell|k}} \to k^*/(k^*)^\nc\right) \to 1.
\end{equation}
By Proposition~\ref{prop:facts-about-l}(iii), the $\nc$-roots  of  unity  in  $k$  and  $\ell$  are  $\pm1$
with the only possible exception when $k = \Q$ and $\ell = \Q[\sqrt{-1}]$.
Therefore, the second group in \eqref{eq:cohom-center-ex-sq} is~$\left\{\pm 1\right\}$
if $\G$ is of type $^2\Dn_{2m+1}$ and is trivial for~$^{3,6}\Dn_4$.

In all the cases we can define the \emph{image} of
$\HH^1(k,\CG)$ in the group $\ell^*/(\ell^*)^\nc$ (for the inner types
we have $\ell=k$). To describe this image
we set $\LL=\ell^*$ in the cases $^2\Dn_{2m}$ and $^1\Dn_{2m+1}$, and
$\LL = \left\{
  x \in \ell^* \; \vert \; \No_{\ell|k}(x) \in (k^*)^\nc \right\}$
in the remaining cases $^2\Dn_{2m+1}$ and $^{3,6}\Dn_4$.
Then, except for the case~$^2\Dn_{2m+1}$, $\LL/(\ell^*)^\nc$
is  an isomorphic  image  of $\HH^1(k,\CG)$.  In the case $^2\Dn_{2m+1}$,
$\HH^1(k,\CG)$ is a $\Z/2\Z$-extension of $\LL/(\ell^*)^\nc$.

\subsection{}
\label{sec:A-instead-whole-H1}
Let $A = \delta(\bG(k))'$ and let
$\AL \le \LL$ be such that the image of $A$ in $\LL/(\ell^*)^\nc$
is equal to $\AL/(\ell^*)^\nc$.
In order to describe the group $A$ we need to solve two problems: the first is related to understanding
the image of $\bG(k)$ and the second is considering the image of $\G(k_{v_0})$
(see Section~\ref{sec:Rohlfs-exct-sequ}).
Both of these images can be described by looking at the archimedean places.

To begin with, let us consider the whole image of
$\delta\!:\!\bG(k)\to\HH^1(k,\CG)$. By (\ref{eq23}) it is equal to
$\ker(\HH^1(k,\CG) \to \HH^1(k,\G))$. The Hasse principle for $\G$
implies that $\HH^1(k,\G)$ is isomorphic
to $\prod_{v\in V_\infty}\HH^1(k_v,\G)$ \cite[Theorem~6.6]{PR}.
Therefore,
\begin{equation}
  \label{eq:delta-bG}
\delta(\bG(k)) = \ker\left( \HH^1(k,\CG) \to {\textstyle\prod_{v\in V_\infty}}\HH^1(k_v,\G)\right).
\end{equation}

Now let us restrict our attention to the group $A$, which means that we have
to add the condition that the elements are in the image of $\G(k_{v_0})$ under the
map $\delta \circ  \phi$. Using~\eqref{eq23} with $K=k_{v_0}$,
we see that these elements are exactly those which are trivial in
$\HH^1(k_{v_0},\CG)$. For the other infinite places $v  \neq v_0$ the
group $\G(k_v) \simeq \Spin(n+1)$ maps
surjectively  onto $\bG(k_v)$  and therefore  $\ker(\HH^1(k_v,\CG) \to
\HH^1(k_v,\G) )= 1$. Thus from~\eqref{eq:delta-bG} we get:
\begin{equation}
  \label{eq:delta-bG-prime}
  A = \ker\left( \HH^1(k,\CG) \to {\textstyle\prod_{v\in V_\infty}}\HH^1(k_v,\CG)\right).
\end{equation}

Let us write $V_\infty = \left\{v_0, v_1 ,\dots\right\}$ for the
infinite  places   of  $k$.   Suppose  that  $\G$   is  not   of
type~$^{3,6}\Dn_4$.  Then the  tensor  product $\ell_{v_i}  = \ell  \otimes
k_{v_i}$ is isomorphic
to $\R$ if $\ell = k$ and isomorphic either to $\C$ or $\R \oplus \R$
according to the type of $\G$ over $k_{v_i}$ when $\ell|k$ is a
quadratic extension. In the latter case,
let us write $\sigma_i$ and $\sigma'_i$ for the two real embeddings of $\ell$ extending
the place $v_i$. In particular, for $x \in \ell$, the notation $x^{\sigma_i} > 0 $ or
$x^{\sigma'_i} > 0$ will make sense.

\begin{prop}
  We have the following description of $A$:
  \label{prop:description-of-A}
  \begin{enumerate}
  \item For  $k \neq  \Q$ we have  $A \simeq \AL/(\ell^*)^\nc$,  and, in
    particular,
    \begin{enumerate}
    \item if $\G$ is of type $^2\Dn_{2m+1}$,
      $\AL = \left\{ x \in \LL \;|\; \emb{\sigma_0}{x} > 0 \right\};$
    \item if $\G$ is of type $^2\Dn_{2m}$,
      $\AL = \left\{ x \in \ell^* \;|\; \emb{\sigma_i}{x},\emb{\sigma'_i}{x} > 0\;
        \forall i \neq 0 \right\}$.
    \end{enumerate}
  \item For $k = \Q$:
    \begin{enumerate}
    \item if $\G$ is of type $^1\Dn_{2m+1}$,
    $\AL = \left\{ x \in \Q^* \;|\; x > 0 \right\}$ and $A \simeq \AL/(\Q^*)^4;$
    \item if $\G$ is of type $^2\Dn_{2m+1}$, $A$ is
      a $\Z/2\Z$-extension  of $\AL/(\ell^*)^4$,  where \\
      $\AL=  \left\{ x \in \LL \;|\; x^{\sigma_0} > 0 \right\};$
    \item if $\G$ is of type $^2\Dn_{2m}$, $\AL = \ell^*$ and $A \simeq \AL/(\ell^*)^2$.
    \end{enumerate}
  \end{enumerate}
\end{prop}
\begin{proof}
  By~\eqref{eq:delta-bG-prime}, an element  of $\HH^1(k,\CG)$ is in $A$
  if and only if its image in $\HH^1(k_v,\CG)$ is trivial for every $v
  \in V_\infty$. The group $\HH^1(k,\CG)$  is described in terms of the
  field $\ell$ and the mapping $\HH^1(k,\CG) \to \HH^1(k_v,\CG)$ can be
  understood when considering the inclusion $\ell \to \ell_v$.
  Let us begin with the type $^2\Dn_{2m+1}$. In this case $\ell_{v_0}
  \simeq \R \oplus \R$ and $\ell_{v_i} \simeq \C$ for $i \neq 0$. The exact
  sequence~\eqref{eq:cohom-center-ex-sq} takes the following form over $v_0$,
  resp. $v_i$:

  \begin{alignat}{4}
     \label{eq:over-v0}
    1 \quad &\to \quad &1 \quad &\to \quad &\HH^1(k_{v_0},\CG) \quad &\to \quad &\left\{ (\R_{>0})^2 , (\R_{<0})^2 \right\} \quad &\to \quad 1;\\
     \label{eq:over-vi}
    1 \quad &\to \quad &\left\{\pm 1\right\} \quad &\to \quad &\HH^1(k_{v_i},\CG) \quad &\to \quad &1 \quad &\to \quad 1.
  \end{alignat}
  Here $(\R_{>0})^2$ and $(\R_{<0})^2$ are meant as the elements
  of $\ell^*_{v_0}/(\ell^*_{v_0})^4 \simeq  (\R^*/\R_{>0})^2$. If $k \neq
  \Q$, there exists a place $v_i \neq v_0$ and hence
  by~\eqref{eq:over-vi}          the         $\Z/2\Z$-extension         is
  killed.  By~\eqref{eq:over-v0},  an element $x \in  \AL$ must respect:
  $x^{\sigma_0}, x^{\sigma'_0} > 0$. But since the norm of $x$ must be
  positive (to be  in $\LL$), $x^{\sigma_0} > 0  $ implies positivity of
  $x^{\sigma'_0}$. These facts together give us {\it (1a)} and {\it (2b)}.

  For the type $^1 \Dn_{2m+1}$ (arising only when $k=\Q$, so
  $V_\infty  = \left\{  v_0  \right\}$ ),  we  have $\HH^1(k_{v_0},\CG)  =
  \HH^1(\R,\mu_4) \simeq  \R^* /  (\R^*)^2$.  An element  $x \in  \Q^*$ has
  then  trivial  image in  $\HH^1(k_{v_0},\CG)$  exactly  when $x$  is
  positive. This proves {\it (2a)}.
  In the case $^2\Dn_{2m}$, the group $\G$ becomes inner over $k_{v_0}$, which means that
  $\ell_{v_0} \simeq  \C$. Then  $\HH^1(k_{v_0},\CG)$ is trivial  and the
  condition at $v_0$ is empty. Hence for $k = \Q$ we
  get { \it (2c)}. For $k \neq \Q$ we need to include the conditions at
  $v_i$ for $i \neq  0$. As $\HH^1(k_{v_i},\CG) \simeq (\R^*/\R_{>0})^2$,
  this gives {\it (1b)}.
\end{proof}

\subsection{}
\label{sec:notations-for-index-comput}

We need to fix some more notations.
Let $T_1$ denote the set of places $v\in V_f$ such that if $\G$ is of
inner type over $k_v$, it does not split, and if it
is of outer type then it is not quasi-split over $k_v$
but splits over its maximal unramified extension $\hat{k}_v$.
In the latter case let $R$ be the subset of $V_f\setminus T_1$
such that $\G$ does not split over $\hat{k}_v$, which can
be also described as the set of places of $k$ which are ramified
in $\ell/k$. We have $T_1 \subset T$, where $T$ is defined
as  in Section~\ref{sec:local-factor} for  some principal  arithmetic
subgroup $\Lambda$ of $\G$.

Let $S$ be a subset of $V_f$. For the subgroup $\AL$ of  $\ell^*$ we denote by
$\AL^S_\nc$ its  subgroup which consists  of the elements $x$  such that
${\tilde  v}(x)\in  \nc\Z$  for every  \emph{normalized}  nonarchimedean
valuation  ${\tilde v}$ of  $\ell$ which  is not  above some  place from
$S$. If $S=\emptyset$ we simply write $\AL_\nc$.
We     introduce     a     refinement     of      the     notation
from Section~\ref{sec:Rohlfs-exct-sequ}:   for the  finite  subset
$S \subset V_f$ we denote by $A_{\xi,S}$ the subgroup of $A$ which acts trivially on
$\Delta_v$  for  every $v  \notin  S$.  Clearly, $A_{\xi,\emptyset}  =
A_\xi$.

\subsection{}
\label{sec:A_xi-from-A_n} By Sections \ref{sec:Rohlfs-exct-sequ}--\ref{sec:Lambda_m}
the index $[\Gamma:\Lambda]$ depends on the order of $A_\xi$.
We can bound this order using results of \cite{BP}.
For  $v \notin S = R \cup T_1$, an element of $A$ acts trivially
on $\Delta_v$ if and only if its image $x$ in $\AL$ has
$\tilde{v}(x) \in \nc \Z$ for every normalized valuation ${\tilde v}$ of $\ell$
extending~$v$ \cite[Lemma~2.3 and Prop.~2.7]{BP}.
This  implies   that  the  image  of  $A_{\xi,S}$  in
$\AL/(\ell^*)^\nc$ is  given by $\AL^S_\nc/(\ell^*)^\nc$.  Therefore, we can
describe the  relation  between   $A_\xi$  and  $\AL_\nc/(\ell^*)^\nc$  by
the following diagram:
\begin{equation*}
  \label{eq:q-q-q}
 \xymatrix{
    *++{A_{\xi,S}}     \ar@{->}[r]^{\bar{q}} &
    *++{\AL^S_{\nc}/(\ell^*)^\nc}  \\
    *++{A_\xi}        \ar@{^{(}->}[u]^{q'}&     *++{  \AL_\nc/(\ell^*)^\nc}
    \ar@{^{(}->}[u]^{q}, \\
  }
\end{equation*}
The integers $q$ and $q'$ stand for the indices of the vertical inclusions.
The order of the kernel of the horizontal map is denoted by $\bar{q}$.
By Proposition~\ref{prop:description-of-A},  $\bar{q}$ is equal
to $1$ with  the only exception of the type $^2\Dn_{2m+1}$
with $k=\Q$, where $\bar{q}=2$.
It is possible  to compute the order of  $\AL_\nc/(\ell^*)$ using the approach
which we will explain in Section~\ref{sec:comput-of-A_n}. The indices  $q$ and  $q'$ are
difficult to evaluate precisely (in particular, for this we would need a counterpart of
\cite[Prop.~2.7]{BP} for $v \in  S$).  However, it is possible
to bound their values:
\begin{prop}
  \label{prop:bound-for-q}
  The indices $q$ and $q'$ satisfy
  \begin{enumerate}
  \item[(1)]
    \begin{enumerate}
    \item $q \le 4^{\# T_1}$, if $\G$ is of absolute type
      $\Dn_{2m+1};$
      \item $q \le 2^{\# R}4^{\# T_1}$, if $\G$ is of type $^2\Dn_{2m}$
      or $^{3,6}\Dn_{4}$.
    \end{enumerate}
  \item[(2)] $q'$ is a factor of $\prod_{v \in S} \# \Xi_v$.
  \end{enumerate}
\end{prop}
\begin{proof}
  The  proof of  the  first assertion  uses  the same arguments as  in
  \cite[Sec. 5]{BP},  where the authors consider the whole $\HH^1(k,\CG)$
  instead of our group $A$. More precisely, for the inner type $^1\Dn_{2m+1}$ we
  obtain  the bound using  the argument of~\cite[Sec.~5.1]{BP}. For  the other
  cases we use the same idea as in the proof of Lemma~5.4 from~\cite{BP}.
  The second assertion is obvious since $A_{\xi,S}$ is
  mapped into $\prod_{v \in S} \Xi_v$ with the kernel $A_\xi$.
\end{proof}

\subsection{}
\label{sec:comput-of-A_n}

We denote by $U_\ell$ the group of units of the ring of integers of $\ell$
and let $U_\AL = \AL \cap U_\ell$. The symbol $\cP_\ell$ stands for the principal
fractional ideal group of $\ell$, $\cC_\ell$ denotes the class  group and its
order  (the  \emph{class  number})  is  denoted by  $h_\ell$.  The  exact
sequence (1) in the proof of \cite[Prop.~0.12]{BP} can be adapted to our setting
by replacing  $\ell^*$ by  $\AL$.  This gives an exact sequence:
$$
1 \to U_\AL \to \AL \to \cP_\AL \to 1,
$$
where $\cP_\AL$ is  a certain subgroup of $\cP_\ell$.  In all cases we
have $(\ell^*)^\nc  \subset \AL$,  so that we  can take  this sequence
modulo $(\ell^*)^\nc$. Then, following the same line of argument as
in the above mentioned proof, we get an exact sequence:
\begin{equation}
  \label{eq:ex-sq-AL}
  1 \to U_\AL/U^\nc_\ell \to \AL_\nc/(\ell^*)^\nc \to \cC_\AL \to 1,
\end{equation}
where  $\cC_\AL$ is a subgroup  of $\cC_\ell$  (it can  be described
explicitly, however, we do not need it here).

\begin{prop}
  \label{prop:bounds-for-index}
  The index of $\Lambda = \Lambda^\mmax$ in its normalizer $\Gamma$
  is bounded by the following value:
  \begin{enumerate}
  \item For $k \neq \Q$:
    \begin{enumerate}
    \item If $\G$ is of type $^2\Dn_{2m+1}$,
      $[\Gamma:\Lambda] \le 2^{d+1} 4^{\#T_1} h_\ell;$
    \item If $\G$ is of type $^2\Dn_{2m}$,
      $[\Gamma:\Lambda] \le  2^{2d-1} 2^{\#R} 4^{\#T_1} h_\ell;$
    \item If $\G$ is of type $^{3,6}\Dn_4$,
      $[\Gamma:\Lambda] \le  2^{3d+1} 2^{\#R} 4^{\#T_1} h_\ell.$
    \end{enumerate}
  \item For $k =  \Q$:
    \begin{enumerate}
    \item If $\G$ is of type $^1\Dn_{2m+1}$,
      $[\Gamma:\Lambda] \le  4^{\#T_1};$
    \item If $\G$ is of type $^2\Dn_{2m+1}$,
      $[\Gamma:\Lambda] \le 8 \cdot 4^{\#T_1} h_\ell;$
    \item If $\G$ is of type $^2\Dn_{2m}$,
      $[\Gamma:\Lambda] \le 4 \cdot 2^{\#R} 4^{\#T_1} h_\ell.$
    \end{enumerate}
  \end{enumerate}
\end{prop}

\begin{proof}
  Let us deal  with the  first case.
  Here we have  $[\Gamma:\Lambda] = \# A_\xi$.
  Proposition~\ref{prop:bound-for-q} and the exact sequence~(\ref{eq:ex-sq-AL})
  allow us to bound the order of $A_\xi$ by $4^{\#T_1}\cdot h_\ell\cdot \#(U_\AL/U^4_\ell)$.
  By Dirichlet's units theorem, $U_\ell$ is a semi-direct
  product of $\mu(\ell) = \left\{\pm 1\right\}$ and $\Z^{d}$
  (where $d+1$ is the number of infinite places
  of $\ell$, as given in Proposition~\ref{prop:facts-about-l}).
  Now we use  $U_\AL \subset U_\LL$ and the fact
  that the norm map $\No_{\ell|k}$ covers $U_k^2$ (taking
  $x$ as a preimage of $x^2$) to get the following bound:
  $$
  \# \left( U_\LL/U_\ell^4 \right) \; \le \;
  \frac{\# \left( U_\ell/U_\ell^4 \right) }{\# \left(  U_k^2/U_k^4 \right)}
  = \frac{2 \cdot 4^{d}}{2^{d-1}}.
  $$
  Finally, restricting our attention to $U_\AL$ instead of $U_\LL$,
  the description of $\AL$ in Proposition~\ref{prop:description-of-A} implies
  that we have to consider only positive elements in $U_\LL$, and so we
  divide the  bound by $2$. This  gives us {\it  (1a)}.  The statement
  {\it (2b)} is obtained by the same  argument with $d=1$,  the only
  difference     comes     from    the value of $q'$ ($2$ instead of
  $1$ when $k=\Q$).

  The remaining cases use only the bound for $q$ in
  Proposition~\ref{prop:bound-for-q}, the bound  for $U_\ell/U^\nc_\ell$
  provided       by      Dirichlet's      units       theorem      and
  Proposition~\ref{prop:facts-about-l}.  Note that for triality forms
  we must  keep in mind that  the index $[\Gamma :  \Lambda]$ is twice
  the order of $A_\xi$.
\end{proof}

\subsection{} Let us point out
that in the proof of Prop.~1.3 \cite[Add.]{B1} it should be written ${\rm Im}(\delta)_\theta$ instead of ${\rm
Im}(\delta)$ everywhere, moreover, in a half of the cases there one has to consider more carefully the set $T$ of places over which the group is non-split. This, however, does not change the result.


\section{\texorpdfstring{Computation of indices $[\Gamma_0:\Lambda_0]$ and $[\Gamma_1:\Lambda_1]$}{Computation of indices}}
\label{sec:index-candidates}
In this section we compute the indices $[\Gamma_0 : \Lambda_0]$ and
$[\Gamma_1 : \Lambda_1]$. The
subgroups          involved          were         defined          in
Section~\ref{section:candidates}.  The  computation  is based  on  the
material of Section~\ref{section:index-estimation}. Similar
kind of arguments have been used in \cite{MS}, where the authors study
groups which act transitively on vertices of Bruhat-Tits buildings.

\subsection{}
\label{sec:index-cand-compact}
Let $k = k_0 = \Q[\sqrt{5}]$ be the field of definition of $\Lambda_0$ and let $\ell = \ell_0$ be
the splitting field attached to $\Lambda_0$.  We have $\D_\ell = \D_{\ell_0} = 275$.  By
construction, the associated group $\G$ is of type $^2\Dn_r$ and we have
$\Lambda_0 = \Lambda^\mmax_0$ which implies that $[\Gamma_0:\Lambda_0] = \# A_\xi$.
By Proposition~\ref{prop:description-of-A}, the group $A$ is isomorphic to
$\AL/(\ell^*)^\nc$ and has two different descriptions
according to the parity of the rank.
Furthermore, as in our case $h_\ell =1 $ the exact sequence~\eqref{eq:ex-sq-AL} shows that
\begin{eqnarray}
  \label{eq:index-in-term-of-q}
  [\Gamma_0 : \Lambda_0] &= & \frac{q}{q'}\; \# \left( U_\AL / U_\ell^\nc
  \right),
\end{eqnarray}
where $q$ and $q'$ are the integers defined in Section~\ref{sec:A_xi-from-A_n}.

The two real embeddings of  $\ell$ correspond to the two $k$-isomorphisms
of $k[\sqrt{\alpha}]$ with $\alpha = 3 + 2 \sqrt{5}$.  These real
embeddings appear over  the place $v_0$ in odd  rank (i.e., in this case $\ell_{v_0} =
\R  \oplus   \R$)  and  over  the   place  $v_1$   in  even  rank
(i.e., here $\ell_{v_1} = \R \oplus \R$).
The field $k[\sqrt{\alpha}]$ has two fundamental units  given by
$$
\tau_1 = \frac{1 + \sqrt{\alpha} }{2} \;\;, \; \;
\tau_2 = \frac{1 - \sqrt{\alpha}}{2}.
$$
We can identify the group $U_\ell/U^\nc_\ell$ with the
representative set
$$
\left\{  \pm \tau_1^i \tau_2^j \; \left| \; 0 \le i,j \le \nc - 1
  \right. \right\}.
$$
The elements which are in  $U_\AL/U_\ell^\nc$ are then
(cf. Proposition~\ref{prop:description-of-A})
\begin{align*}
1 , \tau_1^2 \tau_2^2,  - \tau_1 \tau_2^3 , - \tau_1^3 \tau_2
& \quad \mbox{   when } r = 2m+1; \\
1 ,  - \tau_1 \tau_2
 & \quad \mbox{    when } r = 2m.
\end{align*}
Hence:
\begin{eqnarray}
  \label{eq:number-of-elmt-in-A_n}
\# \left( U_\AL/U_\ell^\nc \right)=
 \left\{
  \begin{array}{cl}
    4 & \mbox{ in case } \Dn_{2m+1},\\
    2 & \mbox{ in case } \Dn_{2m}.
  \end{array}\right.
\end{eqnarray}

\subsection{}
\label{sec:q-q'-for-compact}
We now explain how to determine the values of $q$ and $q'$. We will make
use of the fact stated in \cite[Sec.~2.5]{T2} that for the
ramified place $v_\mathrm{r}$ (defined in the proof of Proposition~\ref{prop:G-compact-quasi-split})
we have $\# \Xi_{v_\mathrm{r}} =  2$, i.e., $\Xi_{v_\mathrm{r}}$ is the whole symmetry group
of the local Dynkin diagram $\Delta_{v_\mathrm{r}}$ of the form
$$\includegraphics{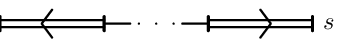}$$
Moreover, we recall that the group $\G_0$ defining  $\Lambda_0$ is such that $T_1  = \emptyset$.

Let us first assume that $r=2m$. In this case $\tau_1 \sqrt{\alpha} \in \AL$ represents
an element of $A_{\xi,R}$ which is clearly not in $\AL_2/(\ell^*)^2$.
This implies that $q = 2$. Moreover, using Hensel's Lemma to
detect congruences modulo square in $\ell_{v_\mathrm{r}}$, we can check that
$-\tau_1 \tau_2$ is another element of $A_{\xi,R}$ which is
mapped to a non-trivial element of $\HH^1(k_{v_\mathrm{r}},\CG) =
\ell_{v_\mathrm{r}}^*/(\ell_{v_\mathrm{r}}^*)^2$.
But since this latter group has order $4$ (see \cite[Ch.~II, Prop.~6]{Lang}),
$A_{\xi,R}$ must map surjectively onto it. In particular, $A_{\xi,R}$ acts
non-trivially on $\Delta_{v_\mathrm{r}}$ and we get $q'=2$.
By Section~\ref{sec:index-cand-compact} this implies that $[\Gamma_0:\Lambda_0] = 2$.

For $r = 2m +1$ by Proposition~\ref{prop:bound-for-q} we have $q=1$.
To see that $A_{\xi,R}$ acts effectively on $\Delta_{v_\mathrm{r}}$ here we cannot use
directly the same argument as for the case $r=2m$. Indeed, $A_{\xi,R}$ cannot be
mapped onto $\HH^1(k_{v_\mathrm{r}},\CG)$ since (contrary to the group $A$) it is
an extension by $\mu_4(k_{v_\mathrm{r}})/N_{\ell|k}(\mu_4(\ell_{v_\mathrm{r}}))
= \left\{\pm 1\right\}$ of the kernel of the norm map
(compare with \eqref{eq:cohom-center-ex-sq}).
But we can check that this subgroup
$\left\{\pm 1\right\} < \HH^1(k_{v_\mathrm{r}},\CG)$ acts trivially
on $\Delta_{v_\mathrm{r}}$. For this we consider the element $(i,i)
\in \mu_4(\overline{k})\times\mu_4(\overline{k}) \simeq
\mathrm{R}_{\ell|k}(\mu_4)(\overline{k})$.
We have $N_{\ell|k}\left( (i,i) \right) = i^2 =
-1$ and hence from the exact sequence \eqref{eq:center-exact-sq}
we see (cf. \cite[Sec. 1.3.2]{PR}) that $-1$ is mapped to the cocycle $a \in \HH^1(k,\CG)$
given by:
\begin{eqnarray*}
	\sigma \in \mathrm{Gal}(\overline{k}|k) &\mapsto& a_\sigma = \left\{
\begin{array}{cl}
	-1 & \mbox{if } \sigma(i) = -i, \\
	1 & \mbox{otherwise.}
\end{array}
\right.
	\label{cocyle-a}
\end{eqnarray*}
We can check that this
cocycle class is represented by $\sigma \mapsto$ $
^{\sigma}\!\tilde{g}\tilde{g}^{-1} \in \CG$, with $\tilde{g} \in
\Spin_{f_0}$ constructed as follows:
\begin{eqnarray*}
\tilde{g} &=& i x y - i x y,
\end{eqnarray*}
for two isotropic vectors $x, y \in (\V,{f_0})$ such that $f_0(x+y)
= 1/2$. We recall here that $\Spin_{f_0}$ is  the subgroup of
the Clifford group which consists of the elements of norm one which give special
orthogonal transformations of $(\V,f_0)$ (see~\cite[Ch.~9, \S3]{Sch}).
A direct computation shows that the image of
$\tilde{g}$ in the group $\SO_{f_0}$ is given by the
matrix $\mathrm{diag}(-1,-1,1,\dots,1)$, with $x,y$ as the first
vectors of the basis completed by orthogonal elements.
Therefore, we see that the image $\phi(\tilde{g})$ (see
Section~\ref{sec:action-H1-on-parahoric}) is in the compact part
of any torus of $\bG(k_{v_\mathrm{r}})$ containing it, which means that
the cocycle $a$ acts trivially on
$\Delta_{v_\mathrm{r}}$ (see \cite[Sec.~2.5]{T2}).  This proves  that $-1 \in
\mu_4(k_{v_\mathrm{r}})$ is trivial in $\Xi_{v_\mathrm{r}}$. Similarly to the case
$r=2m$ we can then check that $A_{\xi,R}$ is mapped onto
$\ker\left(\ell_{v_\mathrm{r}}^*/(\ell_{v_\mathrm{r}}^*)^4 \to
k_{v_\mathrm{r}}^*/(k_{v_\mathrm{r}}^*)^4\right)$.
This shows  that $q' = 2$, so that for all
ranks we have:
\begin{eqnarray}
  \label{eq:index-compact-final}
  [\Gamma_0 : \Lambda_0 ] &=& 2.
\end{eqnarray}

\subsection{}
\label{sec:index-cand-non-compact}

Now consider the group $\Gamma_1$ which is constructed using the group
$\G_1$ defined over $\Q$. As in the cocompact case we have $\Lambda_1 =
\Lambda_1^\mmax$.
Let first $r  = 2m  + 1$.   Then  $\G_1$ is  of inner  type
$^1\Dn_{2m+1}$    and    by    Proposition~\ref{prop:description-of-A}
and Section~\ref{sec:comput-of-A_n}  we  have  $\AL_4/(\Q^*)^4  = 1$.   Now  by
Proposition~\ref{prop:G-local-non-compact}, the set $T_1$ is empty
if $m$ is even, and it follows immediately in this case that $q = q' = 1$.
If  $m$ is odd, $T_1$  consists of a single place  $(2)$.  The elements
$2,4,8 \in  \Q^*$ are three non-trivial elements  of $\AL$ which determine
different   elements    of   $A_{\xi,T_1}$   which    are   not   in
$\AL_4/(\Q^*)^4$.  This allows us to see that $q = 4$. To compute
$q'$ we need to identify $2 \in \Q_2^*/(\Q_2^*)^4$ with an element of
$\bG(\Q_2)$. For this consider a basis $x_1,y_1,\dots,x_r,y_r$ of
$(\V|\Q_2,f_1)$ which satisfies the following conditions:
\begin{itemize}
	\item[(1)] for $j = 1,\dots,r-2$, $x_j$ and $y_j$ are
		isotropic vectors with $f_1(x_j + y_j) = 1/2$;
	\item[(2)] for $j=r-1$ and $j = r$, $x_j$ and $y_j$ are chosen orthonormal.
\end{itemize}
Such a choice of the basis is possible according to the structure
of the space $(\V,f_1)$ over $\Q_2$ (cf. Proposition~\ref{prop:G-local-non-compact}).
With respect to this basis the matrix
\begin{eqnarray*}
	\frac{1}{\sqrt{2}}\left(
	\begin{array}{ccccc}
		\begin{array}{cc}
			2 & 0 \\
			0 & 1
		\end{array}
		& &  &  \\
		 & \ddots &  &  &  \\
		 &  &
		\begin{array}{cc}
			2 & 0 \\
			0 & 1
		\end{array} &  & \\
		 &  &  &
		\begin{array}{cc}
			1 & 1 \\
			-1 & 1
		\end{array}  & \\
		 &  &  & &
		\begin{array}{cc}
			1 & 1 \\
			-1 & 1
		\end{array}
	\end{array}
	\right)
\end{eqnarray*}
gives an element $g$ of $\bG(\Q_2)$ which is contained in the centralizer $Z$ of
a maximal split torus. A preimage $\tilde{g} \in \Spin_{f_1}$ of $g$
under $\phi$ is given by:
\begin{eqnarray*}
	\tilde{g} &=&  \frac{2+\sqrt{2}}{4}  \left(1 + (\sqrt{2}-1) y_{r-1}x_{r-1}\right)
 \left(1 + (\sqrt{2}-1) y_{r}x_{r}\right) \prod_{j=1}^{r-2}
	\left( \sqrt[4]{2} x_j y_j + \frac{1}{\sqrt[4]{2}} y_j
	x_j \right).
\end{eqnarray*}
A direct computation shows that the cocycle $a \in
\HH^1(\Q_2,\CG)$ given by $a_\sigma = {}^\sigma \tilde{g}
\tilde{g}^{-1}$ corresponds exactly to the element $2 \in
\Q_2^*/(\Q_2^*)^4 \simeq \HH^1(\Q_2,\CG)$. We check that $g \not
\in Z_c Z_s$ (see notation in \cite[Sec.~2.5]{T2})
and so $2 \in \HH^1(\Q_2,\CG)$ acts non-trivially on
$\Delta_{(2)}$. Since this local Dynkin diagram has only one
non-trivial symmetry, this implies $q' = 2$ (in case when $r=2m+1$ with
$m$ odd).

For $r  = 2m$  the type  of $\G_1$ is  $^2\Dn_{2m}$ and  the splitting
field $\ell =  \ell_1$ is the field $\Q[\sqrt{-3}]$. We  have $h_\ell = 1$.
By  Proposition~\ref{prop:description-of-A},   $\AL  =  \ell^*$. By
Section~\ref{sec:comput-of-A_n}, we obtain $ \AL_2/(\ell^*)^2  = \left\{
  \pm1  \right\}$.   Again  we  can  compute the value of~$q$:  the  element
$\sqrt{-3}$ represents
a non-trivial element of $A_{\xi,R}$ which is not in
$\AL_2/(\ell^*)^2$, hence $q = 2$. Since $\pm 1, \pm \sqrt{-3}$ are
the four elements composing $\ell_{(3)}^*/\left(\ell_{(3)}^*\right)^2 = \HH^1(\Q_3,\CG)$,
we see that $A_{\xi,R} \to \Xi_{(3)}$ is onto. This implies that
$q' = \#\Xi_{(3)} = 2$.

Similar to Section~\ref{sec:index-cand-compact}, we have
$[\Gamma_1 : \Lambda_1] =  (q/q')\cdot\#\left( U_\AL / U_\ell^\nc\right),$
where now $U_\AL / U_\ell^\nc$ has order $2$ if $r$ is even and is trivial
when $r$ is odd. Hence from the above considerations we get:
\begin{eqnarray}
  \label{eq:index-non-compact-final}
  [\Gamma_1 : \Lambda_1 ] &=&
 \left\{
   \begin{array}{ll}
     1 & \mbox{if } r = 2m+1 \mbox{ with } m \mbox{ even},\\
     2 & \mbox{otherwise}.
   \end{array}\right.
\end{eqnarray}



\section{\texorpdfstring{Uniqueness of $\Gamma_0$ and $\Gamma_1$}{Uniqueness}}
\label{section:uniqueness}
The uniqueness part of Proposition~\ref{prop:G-compact-quasi-split} (resp.
Proposition~\ref{prop:G-local-non-compact}) implies that any
principal arithmetic subgroup of $H$ which has the same global type as
$\Lambda_0$ (resp. $\Lambda_1$) must be commensurable with this
latter group. More precisely, it shows that in each case there exists a uniquely
determined group $\G/k$ and a collection of parahoric subgroups $(\Pa_v)_{v\in V_f}$,
defined up to local conjugations, which are associated to our group $\Lambda_i$.
In this section we will show that this defines $\Lambda_i$ uniquely up to
conjugation by the elements from $\bG(k)$.

\subsection{}
\label{sec:uniqueness-basic}
Let $\cP = (\Pa_v)_{v\in V_f}$ and $\cP' = (\Pa'_v)_{v\in V_f}$ be two coherent collections of parahoric subgroups of $\G$ such that for all $v\in V_f$, $\Pa'_v$ is conjugate to $\Pa_v$ under an element of $\bG(k_v)$. For all but finitely many $v$, $\Pa_v = \Pa_v'$ and they are hyperspecial, hence there exists an element $g\in\bG(\A_f)$ such that $\cP'$ is conjugate to $\cP$ under $g$. Let $\bPa = \prod_{v\in V_f}\bPa_v$ be the stabilizer of $\cP$ in $\bG(\A_f)$. The number of distinct $\bG(k)$-conjugacy classes of coherent collections $\cP'$ as above is the cardinality $c(\bPa)$ of $\mathcal{C}(\bPa) = \bG(k)\backslash \bG(\A_f)/\bPa$, which is called the \emph{class group} of $\bG$ relative to $\bPa$. The \emph{class number} $c(\bPa)$ is known to be finite (see e.g. \cite[Prop.~3.9]{BP}). We need to compute its value for the groups constructed in Section~\ref{section:candidates}.

We first consider a more general setup and then apply the results to our groups.
Recall two isomorphisms (see \cite[Prop.~8.8]{PR}, a minor modification is needed in order to adjust the statement to our setting but the argument remains the same):
\begin{eqnarray*}
\bG(k)\backslash \bG(\A_f)/\bPa & \simeq&  \bG(\A_f)\big/\bG(k)\cdot\bPa;\\
\bG(\A_f)\big/\bG(k)\cdot\bPa & \simeq&  \delta_{\A_f}(\bG(\A_f))
/ \delta_{\A_f}(\bG(k)\cdot \bPa),
\end{eqnarray*}
where $\delta_{\A_f}$ is the restriction of the product map $\prod_v \bG(k_v) \to \prod_v \HH^1(k_v,\CG)$ to the finite ad\`ele group $\bG(\A_f)$.

For every finite place $v$, $\HH^1(k_v,\G)$ is trivial (see
\cite[Theorem~6.4]{PR}) hence $\delta_v: \bG(k_v)\to
\HH^1(k_v,\CG)$ is surjective. Thus the image of
$\delta_{\A_f}(\bG(\A_f))$ can be identified with the restricted
direct product $\prod\nolimits'\HH^1(k_v, \CG)$ with respect to
the subgroups $\delta_v(\bPa_v)$. Also, $\delta_{\A_f}(\bG(k))$
identifies with the image of $\delta(\bG(k))$ under the natural
map $\psi:\HH^1(k,\CG) \to \prod\nolimits'\HH^1(k_v, \CG)$. Hence we have an isomorphism
\begin{eqnarray}
\delta_{\A_f}(\bG(\A_f)) \big/ \delta_{\A_f}(\bG(k)\cdot\bPa)
&\simeq&
{\textstyle\prod'}\:\HH^1(k_v, \CG) \big/
\big((\psi\circ\delta)(\bG(k))\cdot{\textstyle\prod_v}\delta_v(\bPa_v)\big).
\end{eqnarray}

We can summarize the result as follows:

\begin{prop}\label{prop:C(P)} We have $\mathcal{C}(\bPa) \simeq
	\prod'\:\HH^1(k_v, \CG) \big/
	\big((\psi\circ\delta)(\bG(k))\cdot\prod_v\delta_v(\bPa_v)\big).$
\end{prop}

\subsection{}
\label{sec:class-gp-as-ideles}
\newcommand{\Jo}{J_\circ}
To apply this proposition to our groups we will describe
$\prod'\HH^1(k_v,\CG)$ and its quotient using id\`ele groups
(see e.g. \cite[Ch.~VII, \S3]{Lang} for more about id\`eles).
We first consider the case when $\G$ is not of type
$^2\Dn_{2m+1}$. Then $\HH^1(k_v, \CG)$ is canonically isomorphic to
$(\ell\otimes_k k_v)^*/{(\ell\otimes_k k_v)^*}^\nc$.
Now, the group $\delta_v(\bPa_v)$ corresponds to the
kernel $\HH^1(k_v,\CG)_{\xi_v}$ and by \cite[Lemma 2.3 and Prop.
2.7]{BP} (see also Section~\ref{sec:A_xi-from-A_n}), for all but
a finite number of places (namely, $V_f \setminus(T_1 \cup R)$) this kernel
is given by  the elements having  representatives in $\prod_{w|v}\mathfrak{o}_w^* \subset
(\ell\otimes k_v)^*$, where $\mathfrak{o}_w^*$ is the group of
integer units in $\ell_w$. This shows that the
restricted direct product $\prod'\HH^1(k_v,\CG)$ is isomorphic to
$J_f(\ell)/J_f(\ell)^\nc$, where $J_f(\ell)$ denotes   the group
of finite id\`eles of $\ell$.  We shall denote by
$J_{\bPa}/J_f(\ell)^\nc$ its subgroup corresponding to the
product $\prod_v \delta_v(\bPa_v)$ (where $J_{\bPa} \subset
J_f(\ell)$ contains $J_f(\ell)^\nc$).  From the above
explanation it follows that the group
$J_{\bPa}/J_f(\ell)^\nc$ differs only by a finite number of factors from the
group $\Jo(\ell)J_f(\ell)^\nc/J_f(\ell)^\nc$, where
$\Jo(\ell) < J_f(\ell)$ denotes the subgroup of finite id\`eles which are given by units at all
places:
\begin{eqnarray*}
	\Jo(\ell) &=& \prod_{w \in V_f(\ell)} \mathfrak{o}^*_w.
\end{eqnarray*}

We consider now the case $^2\Dn_{2m+1}$ which arises only
for half of the groups $\G_0$. In this case the group
$\HH^1(k_v,\CG)$ was described in Section~\ref{sec:cohom-center}
by the following exact sequence:
$$
1 \to \mu_\nc(k_v) / \No_{\ell|k}(\Res_{\ell|k}(\mu_\nc)(k_v)) \to \HH^1(k_v,\CG) \to
    (\ell\otimes_k k_v)^*/{(\ell\otimes_k k_v)^*}^\nc \stackrel{\No_{\ell|k}}\longrightarrow k_v^*/{k_v^*}^\nc.
$$
We note that by \cite[Sec.~5.3]{BP} together with the argument
from Section~\ref{sec:q-q'-for-compact}, for every $v\in V_f$,
the group $\mu_\nc(k_v) /
\No_{\ell|k}(\Res_{\ell|k}(\mu_\nc)(k_v))$ acts trivially on the
local Dynkin diagram $\Delta_v$ and hence it is contained in
$\delta_v(\bPa_v)$. Therefore, using the same argument as above,  we see
that the quotient $\prod'\HH^1(k_v,\CG) /
\prod_v\delta_v(\bPa_v)$ can be described as a quotient of
$J_f(\LL)/J_f(\ell)^\nc$ by some subgroup $J_\bPa/J_f(\ell)^\nc$
(with $J_f(\ell)^\nc \subset J_\bPa \subset J_f(\LL)$).
Here the id\`ele group  $J_f(\LL)$ is defined by
\begin{eqnarray*}
	J_f(\LL) &=& \left\{ (x_w)_w \in J_f(\ell) \;|\;
	N_{\ell|k}(x_w) \in (k_v^*)^4 \; \mbox{ for all } w|v ,   \;v \in
	V_f\right\},
\end{eqnarray*}
and the group $\LL$ introduced in Section~\ref{sec:cohom-center}
can be seen as a subgroup of $J_f(\LL)$. We can also consider
$\Jo(\LL) = \Jo(\ell) \cap J_f(\LL)$.

Recall that for the cases different from $^2\Dn_{2m+1}$ we
have $\LL = \ell^*$. We will unify the notation for id\`eles in
the rest of the section and denote by $J_f(\LL)$ (resp. $\Jo(\LL)$) the group $J_f(\ell)$
(resp. $\Jo(\ell)$) if the group $\G$ is not of type $^2\Dn_{2m+1}$.

\subsection{}

By Section~\ref{sec:A-instead-whole-H1},
$\delta(\bG(k))\supset A$, with the group $A$
introduced there. By Proposition~\ref{prop:description-of-A}, in
all the cases except when $k = \Q$ and $\G$ is of type
$^2\Dn_{2m+1}$, the group $A$ is a subgroup of
$\ell^*/(\ell^*)^\nc$. By
Proposition~\ref{prop:G-local-non-compact}, when $k = \Q$ and $r
= 2m+1$ our group $\G$ ($=\G_1$) is inner and hence the excluded
case does not occur. So the group $A$ can be identified with
$\AL/(\ell^*)^\nc$, where $\AL$ is given by
Proposition~\ref{prop:description-of-A} {\it (1a), (1b), (2a)} or
{\it (2c)}.

The group $J_\bPa$ introduced in
Section~\ref{sec:class-gp-as-ideles} can differ from $\Jo(\LL)$
at the direct factors corresponding to the places of $\ell$ above the
finite set $R \cup T_1$. But for groups $\G_0$ (resp. $\G_1$)
under consideration, Section~\ref{sec:q-q'-for-compact} (resp.
Section~\ref{sec:index-cand-non-compact}) shows that we have an
inclusion
\begin{eqnarray}
	\AL \Jo(\LL)J_f(\ell)^\nc/J_f(\ell)^\nc  &\subset& \AL
	J_\bPa /J_f(\ell)^\nc.
	\label{AA-cover-J_bPa}
\end{eqnarray}
It follows then from Proposition~\ref{prop:C(P)} and
Section~\ref{sec:class-gp-as-ideles} that the class group
$\mathcal{C}(\bPa)$ is a quotient of the group
\newcommand{\Co}{\mathcal{C}_\circ}
\begin{eqnarray}
\Co &=& (J_f(\LL)/J_f(\LL)^\nc)/(\AL
\Jo(\LL)J_f(\ell)^\nc/J_f(\ell)^\nc) \nonumber \\
	&\simeq& J_f(\LL)/\AL \Jo(\LL) J_f(\ell)^\nc
	\label{eq:Co}.
\end{eqnarray}

\subsection{}
To prove uniqueness it is  sufficient to show that in each case under consideration
the group $\Co$ is trivial. We have the following commutative diagram with an exact
row:
\begin{equation}
	\xymatrix{
	*++{\AL} \ar@{->}[r] \ar@{->}[rd]_{p_1} & *++{J_f(\LL)/\Jo(\LL)J_f(\ell)^\nc}
	\ar@{->}[r] & *++{\Co} \ar@{->}[r] & 1 \\
	& \LL/U_\LL \ar@{->}[u]_{p_2} &
	}
	\label{exact-seqn-Co}
\end{equation}

We need to show the surjectivity of the first horizontal map. To this end
consider the composition $p_2 \circ p_1$. Since $h_\ell =1$ the surjectivity
of $p_2$ can be proved directly, expressing local uniformizers as
the elements in $\ell$.
To prove the surjectivity of $p_1$ we must show that $\AL U_\LL =
\LL$. This is shown if $\LL/\AL$ has the same number of elements
as $U_\LL/U_\AL$. With the description of $\AL$ given in
Proposition~\ref{prop:description-of-A} and the information about
$U_\AL$ from Section~\ref{sec:index-cand-compact} we can easily  check
this equality for all cases. This finishes the proof of uniqueness.

\subsection{}
We showed that the groups $\Lambda_0$ and $\Lambda_1$ are defined uniquely up to conjugation in $\bH$.
It implies that the same is true for their normalizers $\Gamma_0$ and $\Gamma_1$. In the next three sections we will prove that these groups have minimal covolume among the corresponding lattices. The results of this section will then imply the uniqueness of the minimal orbifolds in Theorems 1 and 2.

\section{Proof of Theorem 1: odd rank}
\label{section:proof-compact-odd}
In this section we give a proof of Theorem 1 for the rank $r =
2m+1$. Here we assume that the minimal covolume lattice $\Gamma$
is associated to an algebraic group
$\G$ which is defined  over $k  \neq \Q$  and has the splitting
field   $\ell$.   By Section~\ref{data-from-admissibility}, $\G$ is of type
$^2\Dn_{2m+1}$. Moreover,  we can assume that  $\Gamma = \Gamma^\mmax$
(see Section~\ref{sec:Lambda_m}). We would like to show that $\G = \G_0$.
The results of Sections~\ref{sec:constr-Gamma-0} and \ref{sec:Lambda_m} would
then imply immediately that the minimal covolume is attained on the
lattice $\Gamma_0$.

\subsection{}
\label{sec:T1-minimal}

Let us first consider the case $(k,\ell) = (k_0, \ell_0)$.
Since the rational primes $2$ and $3$ are inert in $k_0$, the
cardinality of the residue field $\fv$ is at least $4$ for
each  $v  \in  V_f(k)$.  By~\eqref{eq:bound-for-lambda-factor},  this
implies  that $\lambda_v \ge 18$
for  every  possible  lambda  factor  appearing  in  the  covolume  of
$\Gamma$.  Then if  $T_1 \neq  \emptyset$ we  see
that  the  covolume  of  $\Gamma$  is strictly  bigger  than  the
covolume of $\Gamma_0$. Indeed, the index $[\Gamma:\Lambda]$ can differ from
$[\Gamma_0:\Lambda_0]$ only by  the factors $q$ and $q'$  (the order of
$\AL_\nc/(\ell^*)^\nc$ depends only on $\ell$), and by
Proposition~\ref{prop:bound-for-q}   we  have   $q   \le  \prod_v
\lambda_v$.
Hence if $(k,\ell) = (k_0,
\ell_0)$, we must  have $T_1 = \emptyset$. By the uniqueness part
of Proposition~\ref{prop:G-local-non-compact} we then have $\G=\G_0$.

Our discussion shows that the proof of the theorem reduces to proving the following
statement: \emph{If $\Gamma $ is
a cocompact arithmetic subgroup of the smallest covolume then the fields $k$, $\ell$
associated to it are the fields $k_0$, $\ell_0$.}
The rest of this section is concerned with
establishing this fact. The idea is to obtain bounds for the discriminants
$\D_k$ and $\D_\ell$ from the inequality $\mu(H/\Gamma) \le
\mu(H/\Gamma_0)$ and the upper bound for $\mu(H/\Gamma_0)$, which we
get from the estimates \eqref{eq:bound_compact} and \eqref{eq:bound_compact-1.2}
(when $r\le  16$)  along with $[\Gamma_0:\Lambda_0]  = 2$.

\subsection{}
\label{sec:bound-for-volume-odd}

From the volume formula and Proposition~\ref{prop:bounds-for-index}
it follows that for $G$ of type $^2\Dn_{2m+1}$:
\begin{eqnarray*}
  \label{eq:bound-for-volume-odd}
  \mu(H/\Gamma) & \ge&\frac{1}{4^{\#T_1} h_\ell 2^{d+1}}
  \D_k^{r^2-r/2} \left( \frac{\D_\ell}{\D_k^2} \right)^{r-1/2}
  C(r)^d \cE(\cP),
\end{eqnarray*}
where $\cP$ is a coherent collection of parahoric subgroups
defining a principal arithmetic subgroup whose normalizer is $\Gamma$.
By~\cite[Appendix C]{BP}, we have $4^{- \#T_1} \cE(\cP) > 1$ (this
value corresponds to $\prod f_v$ in the latter reference), which gives:
\begin{eqnarray}
  \label{eq:bound-for-volume}
  \mu(H/\Gamma) &>& \frac{1}{ h_\ell 2^{d+1}}
  \D_k^{r^2-r/2} \left( \frac{\D_\ell}{\D_k^2} \right)^{r-1/2}
  C(r)^d.
\end{eqnarray}
At this point we need to bound $h_\ell$, and for this we will
use  the  Brauer-Siegel  theorem  similarly  to the  way  it  is  done
in~\cite[proof of Prop.~6.1]{BP}.
Since $\ell$ has signature $(2,d-1)$, the Brauer-Siegel
theorem implies that for any $s > 1$:
\begin{eqnarray}
  \label{eq:Brauer-Siegel}
  h_\ell R_\ell & \le & 2 s (s-1) 2^{-2} \Gamma\left(\frac{s}{2}\right)^2
  \Gamma(s)^{d-1} (2^{-2d+2} \pi^{-2d} \D_\ell)^{s/2} \zeta(s)^{2d},
\end{eqnarray}
where $R_\ell$ is the regulator of $\ell$. In~\cite{F}
it is proved that, with the exception of three fields, the regulator
of a number field is greater than $1/4$. Since these three exceptional fields
are totally complex, they cannot arise as $\ell$. With $R_\ell \ge 1/4$
and taking $s=2$ in~\eqref{eq:Brauer-Siegel}, we get
\begin{eqnarray}
  \label{eq:bound-for-h_l}
  h_\ell & \le & 16 \left( \frac{\pi}{12} \right)^{2d} \D_\ell.
\end{eqnarray}
Combining~\eqref{eq:bound-for-volume}
and~\eqref{eq:bound-for-h_l}, and in the second step using
$\D_\ell \ge \D_k^2$, we can write:
\begin{eqnarray}
  \label{eq:bound-for-volume-without-hl}
  \mu(H/\Gamma) &>& \frac{1}{32}\;
  \D_k^{r^2-5r/2+1}  \D_\ell^{r-3/2} a(r)^d\\
  \label{eq:bound-for-volume-without-hl-without-D_l}
   &\ge& \frac{1}{32}\;
  \D_k^{r^2-r/2-2}  a(r)^d,
\end{eqnarray}
where
$$
a(r) = \frac{1}{2} \left(\frac{12}{\pi} \right)^2 C(r).
$$

\subsection{}
\label{sec:finite-number-of-k}

We now reduce to a finite number of possibilities for
the field $k$ by bounding the discriminant $\D_k$.
For $r \ge 15$ the  factor $a(r)$  is greater than
one, so its power to $d$ can be replaced
in~\eqref{eq:bound-for-volume-without-hl-without-D_l} by a square
without changing the inequality.
Comparing  this value with the upper bound for the covolume
of $\Gamma_0$ we get:
\begin{eqnarray}
  \label{eq:restr-r-bigger-than-15}
  \left( \frac{\D_k}{5 \cdot 11^{1/r}}\right)^{r^2-r/2-2}
  &\le& 128 \left(\frac{\pi}{12} \right)^4
  5^2 \cdot 11^{2/r},
\end{eqnarray}
which is only possible for $\D_k \le 6.1$ if $r\ge15$
(the worst bound appears with the smallest $r$, here $r=15$).
This implies that $k = k_0$.

For $r = 3,5,\dots,13$ the factor $a(r)$ is
less than one and we have to be more careful
in order to exclude fields $k$ of high degrees.
We use \emph{unconditional} (i.e. not depending on the Generalized
Riemann Hypothesis) bounds for
discriminants of totally real number fields given in~\cite[table
4]{Odl}. Taking $d \ge 7$ we have $\D_k > 9.3^d$,
and we can check that for all $r=3,5,\dots, 13$ the factor
$9.3^{r^2-r/2-2} a(r)$ is greater than $1$.
From~\eqref{eq:bound-for-volume-without-hl} we then obtain:
\begin{eqnarray}
  \label{eq:exclude-degree-more-than-7}
  \mu(H/\Gamma) &>& \frac{1}{32}
  \left( 9.3^{r^2-r/2-2}   a(r)\right)^7,
\end{eqnarray}
and we check that for each $r = 3,5,\dots,13$
the right hand side of~\eqref{eq:exclude-degree-more-than-7}
is greater than the covolume of $\Gamma_0$.
So the degree of $k$ must be strictly less than 7.
For the degrees $d=2,\dots,6$ we compare the
bound~\eqref{eq:bound-for-volume-without-hl-without-D_l}
with the covolume of $\Gamma_0$. In each degree
and for each $r$ it gives a bound for $\D_k$.
We list the bounds we get for $r=3$, which are worst than
for other $r$:
\begin{eqnarray*}
  d=2 &:& \D_k \le 22;\\
  d=3 &:& \D_k \le 198;\\
  d=4 &:& \D_k \le 1\;778;\\
  d=5 &:& \D_k \le 15\;956;\\
  d=6 &:& \D_k \le 143\;195.
\end{eqnarray*}
The bounds for $\D_k$ are sufficient to determine from the
tables of totally real number fields (see \cite{tables-bdx} or
\cite{tables-qaos})
precisely which fields may occur as $k$. For $r=7,9,11,15$
we see immediately that the only possibility is $k = k_0$.
For $r=5$ the quadratic field
of discriminant $8$ is not excluded, while for $r=3$ we
are left with fourteen possible fields, all of degree lower
than $5$.

\subsection{}
\label{sec:finite-number-of-l}

We now turn to the problem of excluding all
but a finite number of extensions $\ell|k$.
For $r\ge 15$ and odd we proved that $\D_k =5$
and from~\eqref{eq:bound-for-volume-without-hl} we
can give a bound for $\D_\ell$ similar to the way
we obtained the bound~\eqref{eq:restr-r-bigger-than-15} for $\D_k$.
We obtain that $\D_\ell \le 336$. There are three such
field of signature  $(2,1)$, but only $\ell = \ell_0$  is an extension of
$k_0$.

For $r=3$ the only possible field of degree $d=5$
is $k$ with $\D_k =  14\:641$. Comparing the
bound \eqref{eq:bound-for-volume-without-hl}  for this given  value of
$\D_k$ with $\mu(H/\Gamma_0)$,
we get $\D_\ell \leq 731 \cdot 10^6$. The degree 10 is too large to
find a list of possible number fields in tables from~\cite{tables-bdx}
and~\cite{tables-qaos}. But if we make use of the information
that $\ell$ is an extension of a totally real quintic field,
we can use the results of~\cite{Selm}
where the smallest discriminants of such fields are determined
for each signature. It turns out that our bound for $\D_\ell$ is smaller
than the smallest possible value of the discriminant ($= 1\:332\:031\:009$). This
implies that no extension of $k$ satisfies the inequality, and
hence the quintic field is eliminated.
For $d=4$ we have three  possibilities for $k$ with
$\D_k = 725$, $1\:125$ and $1\:600$. From~\cite{CDM} we
known that here $\D_\ell \ge 4\:286\:875$. Similarly to
what was done for the quintic field above, we
can exclude the case $\D_k = 1\:600$. For the
others two fields this argument does not work. However,
with the method based on the class field theory which is explained
in~\cite{CDM}, we can (using e.g. Pari/GP) find all quadratic extensions
of $k$ with real signature $2$ up to the bound for $\D_\ell$. It turns out that
for $\D_k = 1\:125$ there are no such extensions. For $\D_k = 725$
two extensions $\ell$ are possible (with discriminants
$\D_\ell =  5\:781\:875$ and $9\:986\:875$, resp.), both with class number
$h_\ell = 1$. With this information at hand we can come back to the
inequality~\eqref{eq:bound-for-volume} and refine our bound for $\D_k$. This gives
us the bound $\D_k \le 699$, which excludes the last quartic case
$\D_k = 725$.
For $d=2$, resp. $d=3$, the maximal value for $\D_\ell$
arises when $\D_k = 5$, resp. $\D_k=49$. Using
this value of $\D_k$ in~\eqref{eq:bound-for-volume-without-hl}
and comparing with $\mu(H/\Gamma_0)$ we get the
following bounds for $\D_\ell$:
\begin{eqnarray*}
  d=2 & : & \D_\ell \le 5\:893;\\
  d=3 & : & \D_\ell \le 409\:830.
\end{eqnarray*}
We can  check in the tables from~\cite{tables-bdx}  and~\cite{tables-qaos} that
all fields $\ell$ respecting these bounds have
class  number one.  As  above, this allows us to refine the bound  for $\D_k$  by
putting $h_\ell = 1$ in~\eqref{eq:bound-for-volume}.
We find that the only possible values of $\D_k$ are $5,8,12$ for $d=2$ and
$49$ for $d=3$. For each of these $\D_k$ we
now refine the bound for $\D_\ell$ using $h_\ell=1$.
It turns out that only $k=k_0$ has possible extensions
$\ell$ whose discriminant is small enough, and these possible
extensions are given by $\ell$ with $\D_\ell = 275$, $400$ and $475$.

For $r =5$ the bound for $\D_\ell$ in case $\D_k = 8$ is
small enough to exclude this field. For
$r= 5, \dots, 13$ we then have necessarily $\D_k = 5$,
for which the bounds for $\D_\ell$ allow us to check in a first step that
$h_\ell = 1$. Again this improves the bounds for $\D_\ell$.
From this we find that $\ell$ must be $\ell_0$.

\subsection{}
\label{sec:exclude-last-fields}

To achieve the goal of this section it remains to exclude the cases $\D_\ell
= 400$ or $475$ when $r=3$.
These two discriminants correspond to the extensions $k[\sqrt{\omega}]$ and
$k[\sqrt{\beta}]$ of $k = k_0$, where
$$\omega = \frac{1+\sqrt{5}}{2} \quad \mbox{ and }\quad
\beta=-1 + 2 \sqrt{5}.$$
The fundamental units in  $k[\sqrt{\omega}]$ (resp.  in $k[\sqrt{\beta}]$)
are given by $\tau_1 = \sqrt{\omega}$ and $\tau_2 = 1 + \sqrt{\omega}$
(resp.    $\tau_1    =\frac{1+\sqrt{\beta}}{2}$    and    $\tau_2    =
\frac{1-\sqrt{\beta}}{2}$).
We can compute for these two fields $\ell$ that
$$\#(\AL_4/(\ell^*)^4) = 4.$$
This computation is carried out the same way as it was done for $\ell_0$
in Section~\ref{sec:index-cand-compact}.
This gives  $[\Gamma :  \Lambda] \le
4^{1 + \# T_1}$, and with this new bound we can check that the covolume of
$\Gamma$ is  bigger than  $\mu(H/\Gamma_0)$.

\medskip

We showed that for $r = 2m+1$ the minimal covolume in the cocompact case is
attained on the group $\Gamma_0$ constructed in Section~\ref{section:candidates}.
Together with the index computation in Section~\ref{sec:index-candidates} and
the uniqueness result of Section~\ref{section:uniqueness} it finishes the
proof of Theorem~1 for odd rank.


\section{Proof of Theorem 1: even rank}
\label{section:proof-compact-even}
\subsection{}
\label{sec:general-stuff}

We now deal with the even rank case in Theorem 1. Admissible groups
are of type $^2\Dn_{2m}$ or $^{3,6}\Dn_4$. The proof follows the same
idea as in odd rank and we will allow ourselves to skip some details which
are clear from Section~\ref{section:proof-compact-odd}. Thus, repeating
the argument from Section~\ref{sec:T1-minimal} we  obtain that it is enough to
show that the  minimal covolume subgroups are attached  to $(k,\ell) =
(k_0,\ell_0)$. In particular, in rank $4$ we need to exclude the triality
forms $^{3,6}\Dn_4$. As in Section~\ref{section:proof-compact-odd}, we  denote by
$\Gamma = \Gamma^\mmax$ a subgroup of $\G$ defined over $k \neq \Q$ with
a splitting field $\ell$ and we suppose that $\Gamma$ has minimal
covolume. We will treat the case $^{3,6}\Dn_4$ separately
in Section~\ref{sec:exclude-triality} and assume in the beginning that $\G$ is of
type $^2\Dn_{2m}$.

\subsection{}
\label{sec:bound-for-volume}

The   index
$[\Gamma:\Lambda]$  which appears in  the formula  for $\mu(H/\Gamma)$
is bounded by the integer $2^{2d-1} 2^{\# R} 4^{\#T_1} h_\ell$
(by Proposition~\ref{prop:bounds-for-index}). The  factor $4^{\#T_1}$ can
be eliminated the same way as it was done in Section~\ref{sec:bound-for-volume}. From
\cite[Appendix]{P} we have the inequality
\begin{eqnarray}
  \label{eq:2R-bound}
2^{\# R} & \le &
\frac{\D_\ell}{\D_k^{[\ell:k]}}.
\end{eqnarray}
The bound~\eqref{eq:bound-for-h_l}
for $h_\ell$ is still valid in the current setting. Combining all these
facts together we get
\begin{eqnarray}
  \label{eq:even-bound-for-volume-without-hl}
  \mu(H/\Gamma) &>& \frac{1}{8}\;
  \D_k^{r^2-5r/2+3}  \D_l^{r-5/2} a(r)^d\\
  \label{eq:even-bound-for-volume-without-hl-without-D_l}
   &\ge& \frac{1}{8}\;
  \D_k^{r^2-r/2-2}  a(r)^d,
\end{eqnarray}
where now
$$
a(r) = \frac{1}{4} \left(\frac{12}{\pi} \right)^2 C(r).
$$

\subsection{}
\label{sec:bound-for-Dk-even}
For $r \ge  16$ we can check that  $a(r) \ge 1$ and this  allows us to
bound $\D_k$ by an inequality similar
to~\eqref{eq:restr-r-bigger-than-15}. We  get $\D_k \le  5.9$, proving
that $k=k_0$  if $r \ge  16$. For  $r = 4,  \dots, 14$ we  exclude the
fields $k$ of degree $d \ge 5$ using the same trick as
in Section~\ref{sec:finite-number-of-k} (Odlyzko's bound  employed here is $\D_k
\ge             (6.5)^d$             for            $d             \ge
5$). Comparing~\eqref{eq:even-bound-for-volume-without-hl-without-D_l}
with $\mu(H/\Gamma_0)$ we  get bounds for $\D_k$ depending  on $d$ for
each $r = 4, \dots, 14$. In particular, for $r=4$:
\begin{eqnarray*}
  d=2 &:& \D_k \le 12;\\
  d=3 &:& \D_k \le 62;\\
  d=4 &:& \D_k \le 323.
\end{eqnarray*}
We see  that $k$ can be different from  $k_0$ only
when  $r=4$ (with possible quadratic fields of  discriminants $\D_k= 8,12$
and a cubic field of  discriminant $\D_k =  49$), or when  $r=6$ (with
$\D_k = 8$).

\subsection{}
\label{sec:bound-for-Dl-even}
We now check possibilities for the field $\ell$.
For $r \ge 16$ we obtain the bound $\D_\ell \le 479$.
Besides the field $\ell_0$,  we then have a possibility of $\ell$ having discriminant
$400$ or $475$. These fields have class number $h_\ell =
1$ which allows us to get an inequality
\begin{eqnarray*}
 \left( \frac{\D_\ell}{25 \cdot 11} \right)^{r - 3/2} & \le & 8 \cdot 11.
\end{eqnarray*}
From this inequality we obtain that for  $r \ge 16$, $\D_\ell \le 383$,
which implies that $\ell = \ell_0$.

In rank $r = 4$, for the only possible
cubic field $k$  with $\D_k = 49$ we obtain $\D_\ell \le
16\,809$,   and   no any field of signature   $(4,1)$ satisfies this
condition. Similarly, the cases  $\D_k = 12$ for $r = 4$  and $\D_k = 8$
for  $r  = 6$  can  be also eliminated. Therefore, we obtain the bound $\D_\ell  \le
28\,662$ which is valid in rank  $r=4$ (with the worst case corresponding to $\D_k =5$)
as well as in higher ranks. We
can check in the table from~\cite{tables-qaos} that admissible quartic fields
$\ell$ whose  discriminant respects this bound have  class number $h_\ell
\le 3$. This information can  be used in the inequality
\begin{eqnarray}
  \label{eq:bound-with-hl}
  \mu(H/\Gamma)  & \ge  & \frac{1}{2^{2d-1}  h_\ell} \D_k^{r^2-5r/2+3}
  \D_\ell^{r-3/2} C(r)^d,
\end{eqnarray}
to get  the   refinement $\D_\ell \le  2\,064$. All  fields $\ell$
which satisfy this condition  have the class  number $h_\ell  =  1$. Again,  this
improves  the bound for the discriminant to $\D_\ell \le  1\,330$ when
$r=4$ and  $\D_k = 5$. For  $\D_k = 8$ and the same  rank we have
a better bound $\D_\ell \le 224$ which excludes this case
and proves that $k = k_0$ also in rank $r = 4$.  For the field $\ell$ we now have
the possibilities  $\D_\ell = 275, 400,  475, 775$. For  each of these
fields the  extension  $\ell|k_0$ contains  exactly one  ramified
place, which means that $\# R = 1$. The inequality~\eqref{eq:2R-bound}
can be advantageously replaced by $2^{\#R} = 2$. This gives:
\begin{eqnarray}
  \label{eq:bound-with-one-ramif-place}
  \mu(H/\Gamma)     &      \ge     &     \frac{1}{2^3}     5^{r^2-r/2}
  \left(\frac{\D_\ell}{5^2}\right)^{r-1/2} C(r)^2.
\end{eqnarray}
Comparing the value of this  bound with $\mu(H/\Gamma_0)$ we can check
that $\ell = \ell_0$, with the only possible exception in rank $4$ of
the field $\ell$ with $\D_\ell = 400$. Now for this specific $\ell$ we
can compute that $\#\AL_2/(\ell^*)^2 = 2$  (we use the description  of $\ell$
given in Section~\ref{sec:exclude-last-fields}).   Then   we     have
$[\Gamma:\Lambda] \le 4^{1 + \# T_1}$ which allows us to show that
$\mu(H/\Gamma)$ is  bigger than  $\mu(H/\Gamma_0)$ in the case $\D_\ell =
400$. Thus, in rank $r=4$ we have $\ell = \ell_0$ as well.

\subsection{}
\label{sec:exclude-triality}

To complete the proof of Theorem 1 we need to exclude the arithmetic
$7$-orbifolds defined by triality forms.
So assume now that  $\G$ is  of  type $^{3,6}\Dn_4$.  In this case we have
$[\ell:k] = 3$ and $[\Gamma:\Lambda]  \le 2^{3d+1} 2^{\# R} 4^{\# T_1}
h_\ell$.   Moreover,  the   bound~\eqref{eq:bound-for-h_l}   should  be
replaced by $h_\ell \le 16  (\pi/12)^{3d} \D_\ell$. It follows that
\begin{eqnarray*}
  \label{eq:bound-for-volume-triality}
  \mu(H/\Gamma) &\ge& \frac{1}{32} a^d \D_k^{11},
\end{eqnarray*}
where $a = 1/8 (12/\pi)^3 C(4)$. For $d \ge 5$ we can check using
$\D_k  \ge  6.5^d$ that  $\mu(H/\Gamma)$  is  necessarily bigger  than
$\mu(H/\Gamma_0)$ and for  smaller $d$ we have:
\begin{eqnarray*}
  d=2 &:& \D_k \le 15;\\
  d=3 &:& \D_k \le 86;\\
  d=4 &:& \D_k \le 490.
\end{eqnarray*}
In particular, the totally  real field  $k$ cannot be  of degree  $d =
4$. The inequality
\begin{eqnarray*}
  \mu(H/\Gamma) &\ge& \frac{1}{32} a^d \D_k^{13/2} \D_\ell^{3/2},
\end{eqnarray*}
gives us different bounds for $\D_\ell$ for each of the possible $\D_k$. 
For $\D_k=5$, $d = 2$ and $\D_k = 49$, $d= 3$ we have the worst possible 
bounds:
\begin{eqnarray*}
  d=2 &:& \D_\ell \le 445\:619;\\
  d=3 &:& \D_\ell \le 7.7 \cdot 10^6.
\end{eqnarray*}
We can  check in the table from~\cite{tables-qaos} that   all
fields $\ell$ of degree $6$ with at least one real place (condition~(3)
from Proposition~\ref{prop:facts-about-l}) and with $\D_\ell \le 445\:619$
have the class  number  $h_\ell =  1$. This  information
allows us to improve the bound to $\D_\ell \le 20\,165$, and no any admissible field
$\ell$ has the discriminant small enough. Thus the case $d=2$ is excluded.
For $d=3$ the bound we have is in contradiction with Odlyzko's bound $\D_\ell > (6.1)^9$,
given in \cite{Odl-www}.

\medskip

This finishes the proof of Theorem~1.


\section{Proof of Theorem 2}
\label{section:proof-non-compact}

\subsection{}
\label{sec:case-odd-1}

We begin the proof of Theorem 2  with the easiest case $r =
2m+1$     with     $m$      even.     The     same     argument     as
in Section~\ref{sec:T1-minimal} implies  that it is  enough to show
that a minimal covolume arithmetic subgroup $\Gamma = \Gamma^\mmax$ is
necessarily defined by a group $\G$  of inner type. Thus it suffices to
show that $\D_\ell = 1$. From Proposition~\ref{prop:bounds-for-index},
the bound~\eqref{eq:bound-for-h_l} for $h_\ell$  and the
inequality $4^{-\#T_1}\mathcal{E}(\mathcal{P}) > 1$, we get
\begin{eqnarray*}
  \mu(H/\Gamma)  & \ge  & \frac{1}{8  \cdot 16}  \left( \frac{12}{\pi}
  \right )^2 \D_\ell^{r-3/2}.
\end{eqnarray*}
Since $r  \ge 5$ this gives the bound  $\D_\ell \le 2.27$, proving
that $\ell = \ell_1 = \Q$.

\subsection{}
\label{sec:case-odd-3}

Let us assume now that $r =  2m + 1$ with  $m$ odd. Since  the result of
Theorem  2 for  $r=3$ follows  from~\cite{H}, we  allow
ourselves  to present here the  proof  with the  assumption  that $r  \neq
3$. The proof in the case $r=3$ with our method is also possible but it requires
a bit more effort (see~\cite[\S 15.9]{E-PhD}).

So we have $r \ge 7$. We can check that the product
$$
\zeta(2) \zeta(4) \cdots \zeta(2r -2) \zeta(r)
$$
appearing in $\mu(H/\Lambda_1)$ is bounded above by $2$. Moreover, the lambda
factor $\lambda_{(2)}$ can be bounded by $2^{2r-1}$. Thus
dividing   $\mu(H/\Lambda_1)$  by  the   index  $[\Gamma_1:\Lambda_1]$
from~\eqref{eq:index-non-compact-final},
we obtain:
\begin{eqnarray*}
  \mu(H/\Gamma_1) & < & \frac{4^{r-1/2}}{3} C(r).
\end{eqnarray*}
We  can compare  this  value with  the  general bound  for a  subgroup
$\Gamma   =   \Gamma^\mmax$    obtained   by   our   usual   argument
(Proposition~\ref{prop:bounds-for-index}, bound for $h_\ell$, etc.):
\begin{eqnarray*}
  \mu(H/\Gamma)  & \ge & \frac{1}{8 \cdot 16} \left( \frac{12}{\pi}
  \right)^2 \D_\ell^{r-3/2}.
\end{eqnarray*}
For $r \ge 7$, if we suppose that $\Gamma$ is of the smallest covolume, we
must have $\D_\ell  \le 6.3$. So $\D_\ell = 1$  and $\D_\ell = 5$
are the only possibilities (recall that by
Proposition~\ref{prop:facts-about-l},  $\ell$   must   be totally real).   In
particular, $h_\ell = 1$ and
this gives a  better bound for $\mu(H/\Gamma)$, which is enough to get
$\D_\ell \le 4.7$. The latter is satisfied only for $\ell = \Q$.
Hence the group $\Gamma$  of minimal covolume must be associated to
a group $\G$ of inner type $^1\Dn_{2m+1}$.
Thus $\G$ is inner over each place of $\Q$. Moreover, we know that $G =
\Spin_f$ for some quadratic form $f$ over $\Q$. But since over $\R$ the
Hasse symbol of $f$ must by the admissibility condition be equal to $-1$,
the  Hilbert's   reciprocity  \cite[Ch.  VII]{OM}   implies  the
existence of at least one finite place $v = (p)$ where $\G$ is non-split.
This place $v$ appears in the formula for the covolume of $\Gamma$ with
a non-trivial lambda factor
and we can check using~\eqref{eq:lambda-factor-2} that in case $v \neq (2)$ we
have $\lambda_v /4 \ge \lambda_{(2)}$. Thus considering all other non-trivial
$\lambda_v$, we get $\prod_v  \lambda_v   \cdot   4^{-\#T_1}  \ge
\lambda_{(2)}$.  This implies that the smallest volume is obtained when $\G$
is non-split only at $v = (2)$.  This is exactly the case for $\G
= \G_1$, proving the minimality of the covolume of $\Gamma_1$.

\subsection{}
\label{sec:case-even}

We finally deal with the last case $r= 2m$. Again, even though the proof for $r=4$
is possible, we will assume that $r\ge6$ and refer to \cite{H} for the remaining case.
In  the  covolume  of
$\Gamma_1$ we can bound the Euler product by $2$.  Let us consider the lattice
$\Gamma$ of the minimal covolume. By the admissibility condition, the group $\G$ defining
$\Gamma$ must be of outer type $^2\Dn_{2m}$. Comparing
$\mu(H/\Gamma_1)$ with $\mu(H/\Gamma)$ we have:
\begin{eqnarray}
  \label{eq:bound-for-non-compact-2m}
  3^{r-1/2} & > & \frac{1}{4 \cdot 2^{\#R} h_\ell} \D_\ell^{r-1/2}.
\end{eqnarray}
We use  the  inequality $2^{\#R}  \le
\D_\ell$ to simplify the right hand side.
We  can use the bound~\eqref{eq:bound-for-h_l} which is  valid here
but requires a different argument:  Indeed, this  bound was obtained under
the assumption  that $\ell$ contains only  two roots of  unity. In the
current setting it is no longer the case, since e.g. $\ell  = \Q[i]$ has
$4$ roots  of unity. But  this is the  only exception and  $\Q[i]$ has
the class number $1$ and hence satisfies the bound.
Thus for $r \ge 6$ we obtain $\D_\ell \le 8$. It follows that $\#R = 1$ and
$h_\ell = 1$.  Coming back to~\eqref{eq:bound-for-non-compact-2m} with
the new information at hand, we get $\D_\ell  \le 4$. To exclude the case $\D_\ell
=  4$,  we  compute  the  order  of  $\AL_2/(\ell^*)^2$  for  $\ell  =
\Q[i]$.   By Proposition~\ref{prop:description-of-A},  we
have $\AL = \ell$ and using the
method of Section~\ref{sec:comput-of-A_n}    we    obtain    easily    that
$\ell_2/(\ell^*)^2$ has order $2$.
Compared to the general bound for the index in
Proposition~\ref{prop:bounds-for-index}, this improves  the bound by a
factor two:
$[\Gamma:\Lambda] \le  2 \cdot  2 \cdot 4^{\#T_1}$.  With this new bound
for the index we can
show that $\mu(H/\Gamma)$ is bigger than $\mu(H/\Gamma_1)$. It follows
that the group $\Gamma$  of the smallest  covolume  is necessarily associated to
a group $\G$ with the splitting field $\ell = \ell_1$. We
can see that  the smallest covolume is obtained when no lambda factor  appears in the
volume   formula.  This implies that $\G = \G_1$ and thus  completes
the proof of Theorem~2.


\section{Growth of the minimal volume and proof of Theorem 3}
\label{section:growth}
\subsection{}\label{sec:growth-est}
The formulas from Theorems 1 and 2 allow us to investigate the behavior of the minimal volume as a function of the dimension of the space. As before, we will treat the compact and non-compact cases separately.

By Theorem~1 we have
\begin{eqnarray*}
	\vol(O^n_0) &=&   \frac{5^{r^2-r/2}
     \cdot 11^{r-1/2} \cdot (r-1)!}{2^{2r-1} \pi^r} \; L_{\ell_0|k_0}\!(r) \; \prod_{i=1}^{r-1} \frac{(2i -1)!^2}{(2
     \pi)^{4i}} \zeta_{k_0}(2i),
\end{eqnarray*}
where $r = \frac{n+1}2$. The product $Z = L_{\ell_0|k_0}(r)
\zeta_{k_0}(2) \cdots \zeta_{k_0}(2r-2)$ is minimal for $r =3$,
and by evaluation we get $Z > 1.14$. Thus we have the following
lower bound
\begin{eqnarray}
	\vol(O^n_0) &>& 1.14 \cdot  \frac{5^{r^2-r/2}
     \cdot 11^{r-1/2} \cdot (r-1)!}{2^{2r-1} \pi^r} \; \prod_{i=1}^{r-1} \frac{(2i -1)!^2}{(2
     \pi)^{4i}}.
\end{eqnarray}
It is easy to see that for $r\ge 15$, the product $\prod_{i=1}^{r-1} \frac{(2i -1)!^2}{(2\pi)^{4i}} >1$. This leads to a lower bound
\begin{eqnarray*}\vol(O^n_0) &>& (r-1)! \quad \text{ for } r\ge 15,\end{eqnarray*}
which implies that the minimal volume in Theorem~1 \emph{grows super-exponentially} for high enough $n$.

Similar argument works for the volumes of the non-compact orbifolds from Theorem~2. Here we have
\begin{eqnarray}
\vol(O^n_1) &=&
C_1(r)\prod_{i=1}^{r-1}\frac{(2i -1)!}{(2\pi)^{2i}} \zeta(2i)
> 1.78 \cdot C_1(r)\prod_{i=1}^{r-1}\frac{(2i -1)!}{(2\pi)^{2i}},
\end{eqnarray}
where $r = \frac{n+1}2$ and $C_1(r)$ depends on $r\:\mathrm{(mod\ 4)}$. We can show that for $r\ge20$, the product $\prod_{i=1}^{r-1} \frac{(2i -1)!}{(2\pi)^{2i}} > (2(r-1)-1)!$ and $C_1(r) > 1$ if $r\not\equiv 1\:\mathrm{(mod\ 4)}$, $C_1(r) > \frac{1}{2^{r-2}}$, otherwise. Hence in all the cases
\begin{eqnarray*}\vol(O^n_1) &>& r! \quad \text{ for } r\ge 20. \end{eqnarray*}
	Again we see that the minimal volume grows super-exponentially for high enough $n$.

\subsection{}\label{sec:growth-table}
We can compute the actual values of the minimal volume for small dimensions $n$. Table~2 shows the results of a numerical computation which was done using Pari/GP calculator. It gives \emph{approximate} values of the minimal volume for each of the cases.

\begin{table}[ht]
\label{tab:vol_growth}
$$
\small
\def\arraystretch{1.3}
\begin{array}{|l|*{6}{c}}
\hline
\ n         & 5 & 7 & 9 & 11 & 13 & 15 \quad\ \ \\[3pt]
\hline
\ \vol(O^n_0) & 1.53\cdot10^{-3} & 5.45\cdot10^{-4} & 7.20\cdot10^{-3} & 9.87 & 3.14\cdot10^{6} & 4.57\cdot10^{14} \ \ \dots \\
\ \vol(O^n_1) &  3.65\cdot10^{-4} & 1.89\cdot10^{-6} & 9.30\cdot10^{-11} & 3.52\cdot10^{-11} & 1.20\cdot10^{-12} & 3.52\cdot10^{-14} \quad\ \ \\[3pt]
\hline
\end{array}
$$
\vspace{0.5cm}
$$
\small
\def\arraystretch{1.3}
\begin{array}{*{7}{c}|}
\hline
17 & 19 & 21 & 23 & 25 & 27 & 29\\[3pt]
\hline
\dots\ \ 5.36\cdot10^{25} & 8.40\cdot10^{39} & 2.75\cdot10^{57} & 2.79\cdot10^{78} & 1.27\cdot10^{103} & 3.62\cdot10^{131} & 8.79\cdot10^{163} \\
\quad\ \ 2.07\cdot10^{-18} & 2.69\cdot10^{-14} & 3.65\cdot10^{-12} & 1.13\cdot10^{-10} & 6.63\cdot10^{-13} & 1.25\cdot10^{-2} & 1.98\cdot10^{5} \\[3pt]
\hline
\end{array}
$$
\vspace{0.5cm}
\caption{Approximate volumes of $O^n_0$ and $O^n_1$ for small $n$}
\end{table}

Our computation and estimation show that the minimal volume decrease with $n$ till $n = 7$ (resp. $n = 17$) in the compact case (resp. non-compact case). After this it starts to grow eventually reaching a very fast super-exponential growth which we expect from the estimate in Section~\ref{sec:growth-est}.

\subsection{}\label{sec:proof thm3} We are now ready to prove Theorem~3 from the introduction. The ratio $Q(n) = \vol(O^n_0) / \vol(O^n_1)$ can be bounded from below as follows.

Similar to the previous considerations we can show that for odd $n$,
$$\vol(O^n_1) \le \frac{4^{r-1/2}}{2^{r-1}} 2 \prod_{i=1}^{r-1}\frac{(2i -1)!}{(2\pi)^{2i}} = 2^{r+1}\prod_{i=1}^{r-1}\frac{(2i -1)!}{(2\pi)^{2i}}.$$
It implies
$$Q(n) \ge  \frac{5^{r^2-r/2}\cdot 11^{r-1/2} \cdot (r-1)!}{2^{2r-1} \pi^r \cdot 2^{r+1}} \prod_{i=1}^{r-1}\frac{(2i -1)!}{(2\pi)^{2i}}.$$
Now as in Section~\ref{sec:growth-est}, we have $\prod_{i=1}^{r-1}\frac{(2i -1)!}{(2\pi)^{2i}} > 1$ for $r\ge 15$, which implies that for odd $n\ge 29$, we have $Q(n) > (r-1)!$ where $r = \frac{n+1}{2}$. Together with the computational data from Table~2 this finishes the proof of Theorem~3 for odd $n$. The case of even dimension is obtained by an entirely similar application of the results from \cite{B1}.

\end{document}